\documentclass[namedreferences,draft]{kluwer}

\usepackage{bm}
\usepackage{stmaryrd}

\newlength{\minipagewidth}
\newcommand{\field}[1]{\mathbb{#1}}
\newcommand{\R}{\field{R}}
\newcommand{\E}{\field{E}}

\newcommand{\hx}{\widehat{x}}
\newcommand{\hX}{\widehat{X}}

\newcommand{\bp}{{\boldsymbol{p}}}
\newcommand{\bx}{{\boldsymbol{x}}}

\newcommand{\var}{\mathop{\mathrm{Var}}}
\newcommand{\esp}[2]{\mathbb{E}_{#1} \left[ #2 \right]}
\newcommand{\abs}[1][\cdot]{\arrowvert #1 \arrowvert}
\newcommand{\argmax}{\mathop{\mathrm{argmax}}}
\renewcommand{\ge}{\geq}
\renewcommand{\le}{\leq}
\renewcommand{\textsl}{\emph}

\newdisplay{theorem}{Theorem}
\newdisplay{example}{Example}
\newdisplay{proposition}{Proposition}
\newdisplay{corollary}{Corollary}
\newdisplay{remark}{Remark}
\newdisplay{lemma}{Lemma}

\newcommand{\ignore}[1]{}
\begin{document}

\begin{article}

\begin{opening}

\title{Improved Second-Order Bounds for \\ Prediction with Expert Advice\thanks{
An extended abstract appeared in the \textsl{Proceedings of the 18th Annual Conference
on Learning Theory}, Springer, 2005. The work of all authors was supported in part
by the IST Programme of the European Community, under the PASCAL Network of
Excellence, IST-2002-506778.}}

\author{Nicol\`o \surname{Cesa-Bianchi} \email{cesa-bianchi@dsi.unimi.it}}
\institute{DSI, Universit\`a di Milano, via Comelico 39, 20135 Milano, Italy}
\author{Yishay \surname{Mansour}  \email{mansour@cs.tau.ac.il} \thanks{The
work was done while the author was a fellow in the Institute
of Advance studies, Hebrew University. His work was also supported
by a grant no.\ 1079/04 from the Israel Science Foundation and an IBM faculty
award.}}
\institute{School of computer Science, Tel-Aviv University, Tel Aviv, Israel}
\author{Gilles \surname{Stoltz} \email{gilles.stoltz@ens.fr}}
\institute{D\'epartement de Math\'ematiques et Applications, Ecole Normale Sup\'erieure, 75005 Paris, France}

%NCB I have rewritten this, but I am not happy with the last sentence.
\begin{abstract}
This work studies external regret in sequential prediction games
with both positive and negative payoffs. External regret measures
the difference between the payoff obtained by the forecasting
strategy and the payoff of the best action. In this setting, we
derive new and sharper regret bounds for the well-known
exponentially weighted average forecaster and for a new forecaster
with a different multiplicative update rule. Our analysis has two
main advantages: first, no preliminary knowledge about  the payoff
sequence is needed,
%GS
not even its range;
second, our bounds are expressed in terms of
sums of squared payoffs, replacing larger first-order quantities
appearing in previous bounds. In addition, our most refined bounds
have the natural and desirable property of being stable under
rescalings and general translations of the payoff sequence.
\end{abstract}

\end{opening}

\bibliographystyle{klunamed}

\section{Introduction}
The study of online forecasting strategies in adversarial settings
has received considerable attention in the last few years.
% and not only in the computational learning literature.
One of the goals of the research in this area is the design of
randomized online algorithms that achieve a low external regret;
i.e., algorithms able to minimize the difference between their
expected cumulative payoff and the cumulative payoff achievable
using the single best action (or, equivalently, the single best
strategy in a given class).

If the payoffs are uniformly bounded, and there are finitely many
actions, then there exist simple forecasting strategies whose
external regret per time step vanishes irrespective to the choice of
the payoff sequence. In particular, under the assumption that all
payoffs have the same sign (say positive), the best achieved rates
for the regret are of the order of $\sqrt{X^*}/n$, where $X^*/n$ is
the highest average payoff among all actions after $n$ time steps.
If the payoffs were generated by an independent stochastic process,
however, the tightest rate for the regret with respect to a fixed
action should depend on the \textsl{variance} (rather than the
average) of the observed payoffs for that action.  Proving such a
rate in a fully adversarial setting would be a fundamental result,
and in this paper we propose new forecasting strategies that make a
significant step towards this goal.

Generally speaking, one normally would
expect any performance bound to be maintained under scaling and translation,
since the units of measurement should not make a difference (for example,
predicting the temperature should give similar performances irrespective
to the scale, Celsius, Fahrenheit or Kelvin, on which the temperature is
measured).
However, in many computational settings this does not hold,
for example in many domains there is a considerable difference
between approximating a reward problem or its dual cost problem
(although they have an identical optimal solution).
Most of our bounds also assume no knowledge of the sequence of the
ranges of the payoffs. For this reason it is important for us to
stress that our bounds are stable under rescalings of the
payoff sequence, even in the most general case of payoffs with
arbitrary signs.
The issues of invariance by translations and rescalings,
discussed more in depth in Section~\ref{s:stab}, show
that---in some sense---the bounds introduced in this paper
are more ``fundamental'' than previous results.
In order to describe our results we first set up our model and notations, and
then we review previous related works.

In this paper we consider the following decision-theoretic variant proposed
by~\inlinecite{FS97} of the framework of prediction with expert advice introduced
by~\inlinecite{LW94} and~\inlinecite{Vov98}.
A forecaster repeatedly assigns probabilities to a fixed set of actions.
After each assignment, the actual payoff associated to each action is revealed
and new payoffs are set for the next round. The forecaster's reward on each
round is the average payoff of actions for that round, where the average is
computed according to the forecaster's current probability assignment.
The goal of the forecaster is to achieve, on any sequence of payoffs,
a cumulative reward close to $X^*$, the highest cumulative payoff among
all actions.
We call regret the difference between $X^*$ and the cumulative
reward achieved by the forecaster on the same payoff sequence.

%GS From that point, I had made (2 months ago) several changes, in order to simplify the introduction
%      (since we discuss all known bounds in detail in Section 2).
In Section~\ref{ScTrIssues} we review the previously known bounds on the regret.
The most basic one, obtained via the exponentially weighted average forecaster
of \inlinecite{LW94} and~\inlinecite{Vov98}, bounds the regret by a quantity of the order of $M \sqrt{n \ln N}$,
where $N$ is the number of actions and $M$ is a known upper bound on the magnitude of payoffs.

In the special case of ``one-sided games'', when all payoffs have
the same sign (they are either always nonpositive or always
nonnegative), \inlinecite{FS97} showed that Littlestone and
Warmuth's weighted majority algorithm~\shortcite{LW94} can be used
%NCB we repeat exactly the same sentence in the next section
% as a basic ingredient to construct a forecasting strategy achieving
to obtain
a regret of the order of $\sqrt{M |X^*|\ln N}+M\ln N$.
(If all payoffs are nonpositive, then the absolute value of each payoff is called
\textsl{loss} and $|X^*|$ is the cumulative loss of the best
action.) By a simple rescaling and translation of payoffs, it is
possible to reduce the more general ``signed game'', in which each
payoff might have an arbitrary sign, to either one of the one-sided
games, and thus, bounds can be derived using this reduction. However
the transformation also maps $|X^*|$ to either $Mn + X^*$ or $Mn - X_n^*$,
thus significantly weakening the attractiveness of such a bound.

Recently, \inlinecite{AlNe04} proposed a direct analysis
of the signed game avoiding this reduction.
They proved that weighted majority (used in conjunction with
a doubling trick) achieves the following: on any sequence of payoffs
there exists an action $j$ such that the regret is at most of order
$\sqrt{M(\ln N)\sum_{t=1}^n |x_{j,t}|}$, where $x_{j,t}$ is the
payoff obtained by action $j$ at round $t$, and $M =
\max_{i,t}|x_{i,t}|$ is a known upper bound on the magnitude of
payoffs. Note that this bound does not relate the regret to the sum
$A_n^* = |x_{j^*,1}|+\cdots+|x_{j^*,n}|$ of payoff magnitudes for the
optimal action $j^*$ (i.e., the one achieving $X_n^*$).
%YM Seems like out of context here!
%GS It isn't. The added sentence
%      'Using a standard doubling trick can derive here a bound of $O(\sqrt{M
%      |X^*|\ln N}+M\ln N\ln |X^*|)$'
%      is not true, see the subsection I rewrote in Section 2 about first-order
%      bounds in signed games.
In particular, the bound of order
%GS Corrected a type: it's not |X_n^*| (which is monotone) but A_n^*
$\sqrt{M A_n^* \ln N}+M\ln N$ for one-sided games is only
obtained if an estimate of $A^*_n$ is available in advance.

In this paper we show new regret bounds for signed games. Our analysis has two
main advantages: first, no preliminary knowledge about the payoff magnitude $M$ or about
the best cumulative payoff $X^*$ is needed; second, our bounds are expressed in terms of
sums of squared payoffs, such as $x_{j,1}^2+\cdots+x_{j,n}^2$ and related forms.
These quantities replace the larger terms $M (|x_{j,1}|+\cdots+|x_{j,n}|)$ appearing
in the previous bounds. As an application of our results we obtain, without any
preliminary knowledge on the payoff sequence, an improved regret bound for
one-sided games of the order of $\sqrt{(Mn - |X^*|) (|X^*| / n) (\ln N)}$.

Some of our bounds are achieved using forecasters based on weighted
majority run with a dynamic learning rate. However, we are able to
obtain second-order bounds of a different flavor using a new
forecaster that does not use the exponential probability assignments
of weighted majority. In particular, unlike virtually all previously
known forecasting schemes, the weights of this forecaster cannot be
represented as the gradient of an additive potential (see the
%YM-2-2-06
%monography
monograph
by~\opencite{CBL05} for an introduction to
potential-based forecasters).

\section{An overview of our results}
\label{ScTrIssues}

We classify the existing regret bounds as zero-, first-, and second-order bounds.
%GS I corrected the definitions so that they match our bounds.
%      For instance, unfortunately, not all our first-order bounds depend on $X^*$ (only those
%      which have enough prior information do).
A zero-order regret bound depends only on the number of time steps and on upper bounds on the
individual payoffs. A first-order bound has a main term that depends on a sum of payoffs,
while the main term of a second order bound depends on a sum of squares of the payoffs.
In this section we will also briefly discuss the information which
the algorithms require in order to achieve the bounds.

%GS Put back this paragraph about notation. It's necessary since
%      most notation are not introduced in the introduction anymore (I had moved
%      this paragraph from the introduction to Section 2).
We first introduce some notation
and terminology. Our forecasting game is played in rounds. At each time step $t=1,2,\ldots$
the forecaster computes an assignment $\bp_t = (p_{1,t},\ldots,$ $p_{N,t})$
of probabilities over the $N$ actions. Then the payoff vector
$\bx_t = (x_{1,t},\ldots,x_{N,t}) \in \R^N$ for time $t$ is revealed
and the forecaster's reward is $\hx_t = x_{1,t}p_{1,t} + \cdots + x_{N,t}p_{N,t}$.
We define the cumulative reward of the forecaster by $\hX_n = \hx_1 + \cdots + \hx_n$
and the cumulative payoff of action $i$ by $X_{i,n} = x_{i,1}+\cdots+x_{i,n}$.
For all $n$, let $X_n^* = \max_{i=1,\ldots,N} X_{i,n}$ be the cumulative payoff
of the best action up to time $n$.
The forecaster's goal is to keep the \textsl{regret} $X_n^*-\hX_n$
as small as possible uniformly over $n$.

The one-sided games mentioned in the introduction are the \textsl{loss game}, where
$x_{i,t} \le 0$ for all $i$ and $t$, and the \textsl{gain game}, where
$x_{i,t} \ge 0$ for all $i$ and $t$.
We call \textsl{signed game} the setup in which no assumptions are made on the
sign of the payoffs.

\subsection{Zero-order bounds}
\label{s:zerobd}
We say that a bound is of order zero whenever it only depends on
bounds on the payoffs (or on the payoff ranges)
and on the number of time steps $n$. The basic
version of the exponentially weighted average forecaster of
\inlinecite{LW94} ensures that
the order of magnitude of the regret is
$M \sqrt{n \ln N}$
%GS Note: no extra +M\ln N needed here.
where $M$ is a bound on the payoffs: $\abs[x_{i,t}] \leq M$
for all $t \ge 1$ and $i=1,\ldots,N$.
(Actually, the factor $M$ may be replaced by a
bound $E$ on the \textsl{effective ranges} of the payoffs,
defined by $\abs[x_{i,t} - x_{j,t}] \leq E$ for all $t \ge 1$ and $i,j = 1,\ldots,N$.)
This basic version of this regret bound assumes that we have prior knowledge of both $n$
and $M$ (or $E$).

In the case when $n$ is not known in advance one can use a doubling
trick (that is, restart the algorithm at times $n=2^k$ for
%YM-2-2-06
%some $k$)
 $k\geq \ln N$)
 and achieve a regret bound of the same order, $M \sqrt{n
\ln N}$ (only the constant factor increases). Similarly, if $M$ is
not known in advance, one can restart the algorithm every time the
maximum observed payoff exceeds the current estimate, and take the
double of the old estimate as the new current estimate. Again, this
influences the regret bound by only a constant factor. (The initial
value of the estimate of $M$ can be set to the maximal value in the
first time step, see the techniques used in
Section~\ref{s:refined}.)

A more elegant alternative, rather than the restarting the algorithm from scratch,
is proposed by \inlinecite{AuCeGe02} who consider a time-varying
tuning parameter $\eta_t \sim (1/M)\,\sqrt{(\ln N)/t}$.
They also derive a regret bound of the order of $M\,\sqrt{n\ln N}$ uniformly over
the number $n$ of steps.
%YM What is the dependance here on $M$ ?!
%GS +M\ln N\ln t not needed neither.
Their method can be adapted along the lines of the techniques of Section~\ref{s:E}
to deal with the case when $M$ (or $E$) is also unknown.

%GS
The results for the forecaster of Section~\ref{s:wm} imply a zero-order bound sharper
than $E \sqrt{n \ln M}$. This is presented in Corollary~\ref{CorRange} and basically
replaces $E\sqrt{n}$ by $\sqrt{E_1^2 + \cdots + E_n^2}$, where $E_t$ is the
{\em effective range} of the payoffs at round $t$,
\begin{equation}
\label{defER}
E_t = \max_{i=1,\ldots,N} x_{i,t} - \min_{j=1,\ldots,N} x_{j,t}~.
\end{equation}
%

%NCB Changed title
% \subsection{First-order bounds on the regret (for one-sided games)}
\subsection{One-sided games: first-order regret bounds}
\label{s:firstbd}
We say that a regret bound is first-order whenever its main
term depends on a sum of payoffs. Since the
payoff of any action is at most $Mn$, these bounds are usually
sharper than zero-order bounds. More specifically, they have the
potential of a huge improvement (when, for instance, the payoff of the best action
is much smaller than $Mn$) while they are at most worse by a constant factor
with respect to their zero-order counterparts.

When all payoffs have the same sign \inlinecite{FS97} first showed
that Littlestone and Warmuth's weighted majority
algorithm~\shortcite{LW94} can be used as a basic ingredient to
construct a forecasting strategy achieving a regret of order
$
\sqrt{M |X_n^*|\ln N}+M\ln N
$
where $\abs[X_n^*]$ is
the absolute value of the cumulative payoff of the best action (i.e., the
largest cumulative payoff in a gain game or the smallest cumulative loss in a loss game).

In order to achieve the above regret bound,
the weighted majority algorithm needs prior knowledge of $\abs[X_n^*]$ (or a bound on
it) and of the payoff magnitude $M$. As usual one can overcome this by
a doubling trick. Doubling in this case is slightly more
delicate, and would result in a bound of the order of
$\sqrt{M |X_n^*|\ln N} + M (\ln Mn) \ln N$.
%YM ??? regret bound. (Need to check what we get with simple doubling!)
%GS Done. Also added what follows:
Here again, the techniques of \inlinecite{AuCeGe02} could be adapted along the lines
of the techniques of Section~\ref{s:wm} to get a forecaster that, without restarting and without
previous knowledge of $M$ and $X_n^*$, achieves a regret bounded by
a quantity of the order of $\sqrt{M |X_n^*|\ln N} + M \ln N$.

%NCB Changed title
%\subsection{First-order bounds on the regret (for signed games)}
\subsection{Signed games: first-order regret bounds}
\label{s:firstsigned}
%YM  I cannot follow the logic of keeping this section in!
%GS   Put back this subsection, since we have to cite Allenberg and Neeman.
%        They were the first ones to give sharp first-order bounds in signed games.
%        I however dropped the sentences about translation issues.
%
As mentioned in the introduction, one can translate a signed game
to a one-sided game as follows. Consider a signed game with payoffs
$x_{i,t} \in [-M,M]$. Provided that $M$ is known to the forecaster,
he may use the translation $x'_{i,t} = x_{i,t}+M$ to convert the
signed game into a gain game. For the resulting gain game,
by using the techniques described above,
one can derive a regret bound of the order of
\begin{equation}
\label{EqTr1Sec2}
    \sqrt{(\ln N) \left( Mn + X_n^* \right)} + M\ln N~.
\end{equation}
Similarly, using the translation
$x'_{i,t} = x_{i,t}-M$, we get a loss game, for which one can derive
the similar regret bound
\begin{equation}
\label{EqTr2Sec2}
    \sqrt{(\ln N) \left( Mn - X_n^* \right)} + M\ln N~.
\end{equation}
The main weakness of the transformation is that the bounds~(\ref{EqTr1Sec2})
and~(\ref{EqTr2Sec2}) are essentially zero-order bounds, though
this depends on the precise value of $X_n^*$.
(Note that when $M$ is unknown, or to get tighter bounds,
one may use the translation
$x'_{i,t} = x_{i,t}-\min_{j=1,\ldots,N} x_{j,t}$ from signed games
to gain games, or the translation $x'_{i,t} = x_{i,t}-\max_{j=1,\ldots,N} x_{j,t}$ from
signed games to loss games.)

Recently, \inlinecite{AlNe04} proposed a direct analysis of the
signed game avoiding this reduction. They give a simple algorithm
whose regret is of the order of $ \sqrt{M A_n^* \ln N}+M\ln N $
where $A_n^* = \abs[x_{k_n^*,1}] + \cdots + \abs[x_{k_n^*,n}]$ is
the sum of the absolute values of the payoffs of the best expert
$k_n^*$ for the rounds $1,\ldots,n$.
%YM-2-2-06 : I do not like forward pointers ... and really this is redundant
%(see the notation of Section~\ref{SecNewBdKn}).
Since $A_n^* =
\abs[X_n^*]$ in case of a one-sided game, this is indeed a
generalization to signed games
%NCB
% of the improvement for small payoffs
of Freund and Schapire's first-order bound for one-sided games.
Though Allenberg-Neeman and Neeman need prior knowledge of both $M$ and $A_n^*$
to tune the parameters of the algorithm, a direct extension of their results
along the lines of Section~\ref{SecNewBdKn} gives the first-order bound
% in case $M$, but not $A_n^*$, is known,
\begin{eqnarray}
\nonumber
\lefteqn{
    \sqrt{M (\ln N) \max_{t=1,\ldots,n} A_t^*} + M\ln N
}
\\* &=&
\label{EqTr3Sec2}
    \sqrt{M (\ln N) \max_{t=1,\ldots,n} \sum_{s=1}^t \abs[x_{k_t^*,s}]} + M\ln N
\end{eqnarray}
which holds when only $M$ is known.

\subsection{Second-order bounds on the regret}
\label{s:impact}
%GS I explain below why the previous section is definitely needed.
%
A regret bound is second-order whenever its main term is a function
of a sum of squared payoffs (or on a quantity that is homogeneous in such a sum).
Ideally, they are a function of
\[
    Q_n^*= \sum_{t=1}^n x_{k^*_n,t}^2~.
\]

%YM I would prefer to drop the following paragraph.
%GS I wouldn't: second-order bounds date back to these incomplete information settings,
%      the only thing is that they were never stated or studied per se.
Expressions involving squared payoffs are at the core of many analyses in the
framework of prediction with expert advice, especially in the
presence of limited feedback. (See, for instance, the bandit
problem, studied by~\opencite{AuCeFrSc02}, and more generally
prediction under partial monitoring and the work
of~\opencite{CeLuSt04}, \opencite{CLS04}, \opencite{PiSc01}.)
However, to the best of our knowledge, the bounds presented here are
the first ones to explicitly include second-order information
extracted from the payoff sequence.
%and to be stable under different transformations.

In Section~\ref{s:refined} we give a very simple algorithm whose
regret is of the order of
$
\sqrt{Q_n^* \ln N}+M\ln N
$.
Since $Q_n^* \leq M A_n^*$, this bound improves on the first-order bounds.
Even though our basic algorithm needs prior knowledge of both $M$ and $Q_n^*$
to tune its parameters, we are able to extend it (essentially by using various doubling
tricks) and achieve a bound of the order of
\begin{equation}
\label{EqTr4Sec2}
\sqrt{(\ln N) \max_{t=1,\ldots,n} Q_t^*} + M\ln N
= \sqrt{(\ln N) \max_{t=1,\ldots,n} \sum_{s=1}^t x_{k_t^*,s}^2} + M\ln N
\end{equation}
without using any prior knowledge about $Q_n^*$.
(The extension is not as straightforward as one would expect, since
the quantities $Q_t^*$ are not necessarily monotone over time.)

%GS I finally moved the old
%      \subsection{The impact of extreme values}
%      here. I think it suits well.
Note that this bound is less sensitive to extreme values.
For instance, in case of a loss game (i.e., all payoffs are nonpositive), $Q^*_t \leq M L_t^*$,
where $L_t^*$ is the cumulative loss of the best action up to time $t$.
Therefore, $\max_{s \leq n} Q^*_s \leq M L_n^*$ and the bound~(\ref{EqTr4Sec2})
is at least as good as the family of bounds called ``improvements for small losses''
%YM-2-2-06: We never mention ``improvements for small losses'' there
(or first-order bounds)
presented in Section~\ref{s:firstbd}.
However, it is easy to exhibit examples where the new bound is far
better by considering sequences of outcomes where there are some
``outliers'' among the $x_{i,t}$. These outliers may raise the
maximum $M$ significantly, whereas they have only little impact on
the $\max_{s \leq n} Q^*_s$.

We also analyze the weighted majority algorithm in Section~\ref{s:wm},
and show how exponential weights with a time varying
parameter can be used to derive a regret bound of the order of
$\sqrt{V_n \ln N}+E\ln N$
where $V_n$ is the cumulative variance of the
%NCB
% \emph{online} algorithm
forecaster's rewards
on the given sequence and $E$ is the range of the payoffs. (Again, we
derive first the bound in the case where the payoff range is known,
and then extend it to the case where the payoff range is unknown.)
The above bound is somewhat different from standard regret bounds
because it depends on the predictions of the
%NCB
% online algorithm.
forecaster.
In Sections~\ref{TranslSec} and~\ref{SecAppli} we show how one can use such a
bound to derive regret bounds which only depend on the sequence of payoffs.

\section{A new algorithm for sequential prediction}
\label{s:refined}
We introduce a new forecasting strategy for the signed game.
In Theorem~\ref{th:doubling-refined2}, the main result of this section, we show
that, without any preliminary knowledge of the sequence of payoffs,
the regret of a variant of this strategy is bounded by a quantity
defined in terms of the sums
$Q_{i,n} = x_{i,1}^2+\cdots+x_{i,n}^2$.
Since $Q_{i,n} \le M (|x_{i,1}|+\cdots+|x_{i,n}|)$, such second-order bounds
are generally better than all previously known bounds
%YM-2-2-06
%(see Sections~\ref{s:firstsigned} and~\ref{s:impact}).
(see Section~\ref{ScTrIssues}).

Our basic forecasting strategy, which we call {\tt prod}$(\eta)$,
has an input parameter $\eta > 0$ and maintains a set of $N$ weights.
At time $t=1$ the weights are initialized with $w_{i,1}=1$ for $i=1,\ldots,N$.
At each time $t=1,2,\ldots$, {\tt prod}$(\eta)$ computes the probability
assignment $\bp_t = (p_{1,t},\ldots,p_{N,t})$, where $p_{i,t} = w_{i,t}/W_t$
and $W_t = w_{1,t}+\cdots+w_{N,t}$.
After the payoff vector $\bx_t$ is revealed, the weights are updated using
the rule $w_{i,t+1} = w_{i,t} (1 + \eta x_{i,t})$.
The following simple fact plays a key role in our analysis.
\begin{lemma}
\label{first}
For all $z \ge -1/2$, $\ln(1+z) \ge z - z^2$.
\end{lemma}
\begin{pf}
Let $f(z) = \ln(1+z) - z + z^2$. Note that
\[
f'(z) = \frac{1}{1+z} - 1 + 2z = \frac{z(1+2z)}{1+z}
\]
so that $f'(z) \leq 0$ for $-1/2 \leq z \leq 0$ and $f'(z) \geq 0$ for $z \geq 0$.
Hence the minimum of $f$ is achieved in $0$ and equals $0$, concluding the proof.
\qed
\end{pf}
We are now ready to state a lower bound on the cumulative reward of {\tt prod}$(\eta)$
in terms of the quantities $Q_{k,n}$.
\begin{lemma}
\label{lm:refined}
Assume there exists $M > 0$ such that the payoffs satisfy $x_{i,t} \ge -M$
for $t=1,\ldots,n$ and $i=1,\ldots,N$.
For any sequence of payoffs, for any action $k$, for any $\eta \leq {1}/({2M})$, and
for any $n \ge 1$, the cumulative reward of ${\tt prod}(\eta)$ is lower bounded as
\[
    \hX_n \ge X_{k,n} - \frac{\ln N}{\eta} - \eta\,Q_{k,n}~.
\]
\end{lemma}
\begin{pf}
For any $k=1,\ldots,N$, note that $x_{k,t} \ge -M$ and $\eta \le 1/(2M)$
imply $\eta x_{k,t} \ge -1/2$. Hence, we can apply Lemma~\ref{first} to
$\eta x_{k,t}$ and get
\begin{eqnarray}
\nonumber \lefteqn{   \ln\frac{W_{n+1}}{W_1}} \\
\nonumber
& & \geq \ln
\frac{w_{k,n+1}}{W_1}
 =  - \ln N + \ln\prod_{t=1}^n (1 + \eta x_{k,t})
= - \ln N + \sum_{t=1}^n \ln(1 + \eta x_{k,t}) \\
& &  \ge
\label{new-lower}
    - \ln N + \sum_{t=1}^n \left(\eta x_{k,t} - \eta^2 x_{k,t}^2\right)
=
    -\ln N + \eta X_{k,n} - \eta^2 Q_{k,n}~.
\end{eqnarray}
On the other hand,
\begin{eqnarray}
\nonumber
     \ln\frac{W_{n+1}}{W_1}
& = &
    \sum_{t=1}^n\ln\frac{W_{t+1}}{W_t}
=
    \sum_{t=1}^n \ln \left( \sum_{i=1}^N p_{i,t}\,(1+\eta x_{i,t}) \right) \\
\label{new-upper}
& = &
    \sum_{t=1}^n \ln\left(1+\eta\sum_{i=1}^N x_{i,t}p_{i,t}\right)
\le
    \eta\hX_n
\end{eqnarray}
where in the last step we used $\ln(1+z_t) \le z_t$ for all $z_t =
\eta\sum_{i=1}^N x_{i,t}p_{i,t} \geq -1/2 $.
Combining~(\ref{new-lower}) and~(\ref{new-upper}), and dividing by
$\eta > 0$, we get
\[
    \hX_n \ge - \frac{\ln N}{\eta} + X_{k,n} - \eta\,Q_{k,n}~.
\]
Our choice of $\eta$ gives the claimed bound.
%GS Yishay, did you write that here intentionally?
%      'Since $\eta \leq 1/2M$, we have that $z_t \geq -M/2M=-1/2$'
%      We already have it 10 lines above.
\qed
\end{pf}
By choosing $\eta$ appropriately, we can optimize the bound as follows.
\begin{theorem}
\label{th:refined}
Assume there exists $M > 0$ such that the payoffs satisfy $x_{i,t} \ge -M$
for $t=1,\ldots,n$ and $i=1,\ldots,N$.
For any $Q > 0$,
if ${\tt prod}(\eta)$ is run with
\begin{equation}
\label{etaTH1}
    \eta = \min\left\{ {1}/{(2M)} , \sqrt{{(\ln N)}/{Q}} \right\}
\end{equation}
then for any sequence of payoffs, for any action $k$, and for any $n \ge 1$
such that $Q_{k,n} \le Q$,
\[
    \hX_n \ge X_{k,n} - \max\left\{ 2\sqrt{Q \ln N}~, 4M\ln N \right\}~.
\]
\end{theorem}

\subsection{Unknown bound on quadratic variation (Q)}
\label{SecNewBdKn}
To achieve the bound stated in Theorem~\ref{th:refined}, the
parameter $\eta$ must be tuned using preliminary knowledge of a
lower bound on the payoffs and an upper bound on the quantities
$Q_{k,n}$. In this and the following sections we remove these
requirements one by one. We start by introducing a new algorithm
that, using a doubling trick over ${\tt prod}$, avoids any
preliminary knowledge of an upper bound on the $Q_{k,n}$.

Let $k_t^*$ be the index of the best action up to time $t$; that is, $k_t^* \in \argmax_k X_{k,t}$
(ties are broken by choosing the action $k$ with minimal associated $Q_{k,t}$).
We denote the associated quadratic penalty by
\[
Q^*_t = Q^*_{k_t^*} = \sum_{s=1}^t x_{k_t^*,s}^2~.
\]
Ideally, our regret bound should depend on $Q^*_n$ and be of the
form ${\sqrt{Q^*_n \ln N}} + M\ln N$. However, note that the
sequence $Q^*_1,Q^*_2,\ldots$ is not necessarily monotone,
%YM
since if at time $t+1$ the best action changes, then $Q_t^*$ and $Q_{t+1}^*$ are
not related.
%as $Q_t^*$ and $Q_{t+1}^*$ cannot be possibly related when the actions achieving the
%largest cumulative payoffs at rounds $t$ and $t+1$ are different.
Therefore, we cannot
use a straightforward doubling trick, as this only applies to monotone sequences.
Our solution is to express the bound in terms of the smallest nondecreasing sequence that
upper bounds the original sequence $(Q_t^*)_{t \geq 1}$.
This is a general trick to handle situations where the penalty terms are not monotone.

Let {\tt prod-Q}$(M)$ be the prediction algorithm that receives a
quantity $M > 0$ as input parameter and repeatedly runs ${\tt prod}(\eta_r)$,
where $\eta_r$ is defined below.
The parameter $M$ is a bound on the payoffs, such that
for all $i=1,\ldots,N$ and $t=1,\ldots,n$,
we have $\abs[x_{i,t}] \leq M$.
% Note that we need the absolute values to be sure that the squared terms Q_{k,n} are
% less than M^2 n (used in the last line of the proof of Th.2) -- In Th.1, we did not need these absolute values.
The $r$-th parameter $\eta_r$ corresponds to the parameter $\eta$ defined in~(\ref{etaTH1})
for $M$ and $Q = 4^r M^2$. Namely, we choose
\[
\eta_r = \min\left\{ {1}/{(2M)} , \sqrt{\ln N}/(2^r M) \right\}~.
\]
We call epoch $r$, $r=0,1,\ldots$, the
sequence of time steps when $\mbox{\tt prod-Q}$ is running ${\tt prod}(\eta_r)$.
The last step of epoch $r \ge 0$ is the time step $t=t_r$
when $Q^*_t > 4^r M^2$ happens for the first time.
When a new epoch $r+1$ begins, ${\tt prod}$ is restarted with parameter
$\eta_{r+1}$.
\begin{theorem}
\label{th:doubling-refined}
Given $M > 0$, for all $n \geq 1$ and all sequences of payoffs bounded by $M$,
i.e., $\max_{1\le i \le N}\max_{1 \le t \le n} |x_{i,t}| \leq M$,
the cumulative reward of algorithm {\tt prod-Q}$(M)$ satisfies
\[
\begin{array}{lllll}
    \hX_n &
\ge &
    X_n^* & - & 8 \sqrt{ (\ln N) \max_{s \leq n} Q^*_s} \\*
    & & & - & 2\, M \Bigl( 1 + \log_4 n + 2 \bigl( 1 + \lfloor (\log_2 \ln N)/2 \rfloor \bigr) \ln N
    \Bigr)~\\
    & = & X_n^* & - & O \left( \sqrt{ (\ln N) \max_{s \leq n} Q^*_s} + M\ln n
    + M \ln N \ln \ln N \right)~.
\end{array}
\]
\end{theorem}
\begin{pf}
We denote by $R$ the index of the last epoch and let $t_R = n$.
If we have only one epoch, then the theorem follows from Theorem~\ref{th:refined}
applied with a bound of $Q = M^2$ on the squared payoffs of the best expert.
Therefore, for the rest of the proof we assume $R \ge 1$.
Let
\[
    X_k^{(r)} = \sum_{s=t_{r-1}+1}^{t_r-1} x_{k,s}~,
\quad
    Q_k^{(r)} = \sum_{s=t_{r-1}+1}^{t_r-1} x_{k,s}^2~,
\quad
    \hX^{(r)} = \sum_{s=t_{r-1}+1}^{t_r-1} \hx_{s}
\]
where the sums are over all the time steps $s$ in epoch $r$ except the last one, $t_r$.
(Here $t_{-1}$ is conventionally set to $0$.)
We also denote $k_r = k^*_{t_r - 1}$ the index of the best overall expert up to time $t_r-1$
(one time step before the end of epoch $r$).
We have that $Q_{k_r}^{(r)} \le Q_{k_r,t_r - 1} = Q^*_{t_r - 1}$.
Now, by definition of the algorithm, $Q^*_{t_r - 1} \leq 4^r M^2$.
Theorem~\ref{th:refined} (applied to time steps $t_{r-1}+1,\ldots,t_r-1$) shows that
\[
    \hX^{(r)} \geq X_{k_r}^{(r)} - \max\left\{ 2\sqrt{4^r M^2 \ln N}~, 4M\ln N \right\}~.
\]
The maximum in the right-hand side equals $2^{r+1} M \sqrt{\ln N}$
when $ r > r_0 = 1 + \lfloor (\log_2 \ln N)/2 \rfloor$.
Summing over $r=0,\ldots,R$ we get
\begin{eqnarray}
\nonumber
\hX_n & = & \sum_{r=0}^R \left(\hX^{(r)} + \hx_{k_r,t_r}\right) \\
\nonumber
& \geq &
\sum_{r=0}^R \left( \hx_{k_r,t_r} + X_{k_r}^{(r)} \right) - 4 ( 1 + r_0)M  \ln N - \sum_{r=r_0+1}^R 2\sqrt{4^r M^2 \ln N} \\
% We can have a small gain of $2^{r_0}$  but it seems negligible.
\nonumber & \geq &
\sum_{r=0}^R \left( \hx_{k_r,t_r} + X_{k_r}^{(r)} \right) - 4 ( 1 + r_0) M \ln N - 2^{R+2} M \sqrt{\ln N} \\
\label{EqCumDT1} & \geq & \sum_{r=0}^R X_{k_r}^{(r)} - (R+1)M - 4
( 1 + r_0) M \ln N - 2^{R+2} M \sqrt{\ln N}~.
\end{eqnarray}
Now, since $k_0$ is the index of the expert with largest payoff up to time $t_0 - 1$,
we have that
$X_{k_1, t_1 - 1} = X_{k_1}^{(0)} + x_{k_1,t_0} + X_{k_1}^{(1)} \leq X_{k_0}^{(0)} + X_{k_1}^{(1)} + M$.
By a simple induction, we in fact get
\begin{equation}
\label{EqCumDT2}
X_{k_R,t_R-1} \leq \sum_{r=0}^{R-1} \left( X_{k_r}^{(r)} + M \right) + X_{k_R}^{(R)}~.
\end{equation}
As, in addition, $X_{k_R,t_R-1} = X_{k^*_{n-1},n-1}$ and $X_{k^*_n,n}$ may only differ by at most $M$,
combining~(\ref{EqCumDT1}) and~(\ref{EqCumDT2}) we have indeed proven that
\[
\hX_n \geq X_{k^*_n,n} - \left( 2 (R+1) M + 4 M ( 1 + r_0) \ln N + 2^{R+2} M \sqrt{\ln N} \right)~.
\]
The proof is concluded by noting first, that $R \leq \log_4 n$, and second that,
as $R \geq 1$, $\max_{s \leq n} Q^*_s \geq 4^{R-1} M^2$ by definition of the algorithm.
\qed
\end{pf}

\subsection{Unknown bound on payoffs (M)}
\label{SecUnkM}

In this section we show how one can overcome the case when there is no a priori
bound on the payoffs.
In the next section we combine the techniques of this section and Section~\ref{SecNewBdKn}
to deal with the case when both parameters are unknown

Let {\tt prod-M}$(Q)$ be the prediction algorithm that receives a
number $Q > 0$ as input parameter and repeatedly runs ${\tt prod}(\eta_r)$,
where the $\eta_r$, $r=0,1,\ldots$, are defined below.
We call epoch $r$ the
sequence of time steps when $\mbox{\tt prod-M}$ is running ${\tt prod}(\eta_r)$.
At the beginning, $r=0$ and $\mbox{\tt prod-M}(Q)$ runs ${\tt prod}(\eta_0)$,
where
\[
M_0 = \sqrt{{Q}/{(4\ln N)}} \qquad \mbox{and} \ \
\eta_0 = {1}/{(2 M_0)} = \sqrt{{(\ln N)}/{Q}}~.
\]
For all $t \geq 1$, we denote
\[
M_t = \max_{s=1,\ldots,t} \max_{i = 1,\ldots,N} 2^{\lceil \log_2 |x_{i,s}| \rceil}~.
\]
The last step of epoch $r \ge 0$ is the time step $t=t_r$
when $M_t > M_{t_{r-1}}$ happens for the first time
(conventionally, we set $M_{t_{-1}}=M_0$).
When a new epoch $r+1$ begins, ${\tt prod}$ is restarted with parameter
$\eta_{r+1} = 1/(2 M_{t_r})$.

Note that $\eta_0 = 1/(2 M_0)$ in round $0$ and $\eta_r = 1/(2 M_{t_{r-1}})$
in any round $r \ge 1$, where
$M_{t_0} > M_0$ and $M_{t_r} \geq 2 M_{t_{r-1}}$ for each $r \geq 1$.
\begin{theorem}
\label{guess-M}
For any sequence of payoffs, for any action $k$, and for any $n \ge 1$
such that $Q_{k,n} \le Q$, the cumulative reward of algorithm {\tt prod-M}$(Q)$
is lower bounded as
\[
% In the case R = 0, we get exactly 2 \sqrt{Q \ln N} from Th.1
    \hX_n \ge X_{k,n} - 2 \sqrt{Q \ln N} - 12 \,M (1+\ln N)
\]
where
$M = \max_{1\le i \le N}\max_{1 \le t \le n} |x_{i,t}|$.
\end{theorem}
\begin{pf}
As in the proof of Theorem~\ref{th:doubling-refined}, we denote by
$R$ the index of the last epoch and let $t_R = n$. We assume $R \ge
1$ (otherwise, the theorem follows directly from
Theorem~\ref{th:refined} applied with a lower bound of $- M_0$ on
the payoffs).
% What follows is needed for later purposes -- In the COLT version, we did not pay attention to the fact
% that sometimes t_R = n just because the last epoch really ends at n in the sense
% that the bound doubles at the last time step
Note that at time $n$ we have either $M_n \le M_{t_{R-1}}$, implying $M_n = M_{t_R} = M_{t_{R-1}}$,
or $M_n > M_{t_{R-1}}$, implying $M_n = M_{t_R} = 2 M_{t_{R-1}}$.
In both cases, $M_{t_R} \ge M_{t_{R-1}}$. In addition, since $R \ge 1$, we also have
$M_{t_{R}} \le 2M$.

Similarly to the proof of Theorem~\ref{th:doubling-refined}, for all epochs $r$ and actions $k$
introduce
\[
    X_k^{(r)} = \sum_{s=t_{r-1}+1}^{t_r-1} x_{k,s}~,
\quad
    Q_k^{(r)} = \sum_{s=t_{r-1}+1}^{t_r-1} x_{k,s}^2~,
\quad
    \hX^{(r)} = \sum_{s=t_{r-1}+1}^{t_r-1} \hx_{s}
\]
where, as before, we set $t_{-1}=0$. Applying
Lemma~\ref{lm:refined}
to each epoch $r = 0,\ldots,R$ we get that
$\hX_n - X_{k,n}$ is equal to
\begin{eqnarray*}
\hX_n-X_{k,n} &= &
\sum_{r=0}^{R} \left(\hX^{(r)} - X^{(r)}_k \right) + \sum_{r=0}^{R} \left( \hx_{t_r} - x_{k,t_r} \right)
\\
& & \ge
    - \sum_{r=0}^R \frac{\ln N}{\eta_r}
    - \sum_{r=0}^R \eta_r Q^{(r)}_{k}
   + \sum_{r=0}^{R} \left( \hx_{t_r} - x_{k,t_r} \right)~.
\end{eqnarray*}
We bound each sum separately.
For the first sum,
since $M_{t_s} \ge 2^{s-r} M_{t_r}$ for each
$0 \leq r \leq s \leq R-1$, we have for $s \leq R-1$,
\begin{equation}
\label{sum2power}
    \sum_{r=0}^s M_{t_r}
\le
    \sum_{r=0}^s 2^{r-s} M_{t_s}
\leq
    2 M_{t_s}~.
\end{equation}
Thus,
\[
    \sum_{r=0}^R \frac{\ln N}{\eta_r}
=
    \sum_{r=0}^R 2 M_{t_{r-1}} {\ln N}
\le
    2\Bigl(M_{t_{-1}} + 2 M_{t_{R-1}}\Bigr)\ln N
\le
    6 M_{t_{R}} \ln N
\]
where we used~(\ref{sum2power}) and $M_{t_{-1}} = M_0 \le
M_{t_{R-1}} \le M_{t_{R}}$. For the second sum, using the fact that
$\eta_r$ decreases with $r$, we have
\[
\sum_{r=0}^R \eta_r Q^{(r)}_{k} \leq \eta_0 \sum_{r=0}^R  Q^{(r)}_{k}
% Note that we do not count x_{t_r}^2
\leq \eta_0 Q_{k,n} \leq \sqrt{\frac{\ln N}{Q}} \, Q = \sqrt{Q \ln N}~.
\]
Finally, using~(\ref{sum2power}) again,
\[
    \sum_{r=0}^R \left| \hx_{t_r} - x_{k,t_r} \right|
\le
    \sum_{r=0}^R 2\,M_{t_r} \leq 2 \left( 2\,M_{t_{R-1}} + M_{t_R} \right) \leq 6\,M_{t_R}~.
\]
The resulting lower bound $\displaystyle{6 M_{t_{R}} (1+\ln N) + \sqrt{Q \ln N}}$
implies the one stated in the theorem by recalling that, when $R \geq 1$, $M_{t_R} \leq 2\,M$.
\qed
\end{pf}

\subsection{Unknown bounds on both payoffs (M) and quadratic variation (Q)}
\label{s:MQ}
%
% In the COLT version of this paper, we had a $\max \{ 1, ... \}$ because we started the doubling trick on $Q_t^*$
% at 1, no matter what $M$ was!
% If we used $4^r Q^*_1$ to do the doubling trick (one of the attempts), then the main root term would be OK, but we would get
% a term of the order of $M \ln (Q^*_n/Q_1^*)$, that is, a term very similar to $\ln V_n/V_1$ we had to get rid of in the analysis of
% the weighted majority algorithm!
%
We now show a regret bound for the case when $M$ and the $Q_{k,n}$ are both unknown.
We consider again the notation of the beginning of Section~\ref{SecNewBdKn}.
The quantities of interest for the doubling trick of Section~\ref{SecNewBdKn} were the
homogeneous quantities $(1/M^2) \max_{s \leq t} Q^*_s$. Here we assume no knowledge of $M$.
We propose a doubling trick on the only homogeneous quantities we have access to, that is,
$\max_{s \leq t} \, (Q^*_s/M_s^2)$, where $M_t$ is defined in Section~\ref{SecUnkM}
and the maximum is needed for the same reasons of monotonicity as in Section~\ref{SecNewBdKn}.

We define the new (parameterless) prediction algorithm {\tt prod-MQ}.
Intuitively, the algorithm can be thought as running, at the low level,
the algorithm $\mbox{\tt prod-Q}(M_t)$.
When the value of $M_t$ changes, we restart $\mbox{\tt prod-Q}(M_t)$, with the new value
% We do not completely restart the algorithm: we do not forget the value of $Q_t^*$
but keep track of $Q_t^*$.
% I would't write what follows, because the algorithm also takes fresh starts within the prod-Q(M_t) steps.
% Or, alternatively, I would upper bounds the total number of fresh starts in the nested doubling.
% '(Note that the number of restarts is at most $\ln M$, and each time the value of $M$ at least
% doubles.)'

Formally, we define the prediction algorithm {\tt prod-MQ} in the
following way. Epochs are indexed by pairs $(r,s)$. At the beginning
of each epoch $(r,s)$, the algorithm takes a fresh start and runs
${\tt prod}(\eta_{(r,s)})$, where $\eta_{(r,s)}$,
for $r=0,1,\ldots$ and $s=0,1,\ldots$, is defined by
\[
\eta_{(r,s)} = \min\left\{ {1}\Big/{\Bigl(2M^{(r)}\Bigr)},\,
\sqrt{\ln N}\Big/\Bigl(2^{S_{r-1}+s} M^{(r)}\Bigr) \right\}
\]
and $M^{(r)},\, S_r$ are defined below.

At the beginning, $r=0$, $s=0$, and since \texttt{prod}$(\eta)$ always sets $\bp_1$
to be the uniform distribution irrespective to the choice of $\eta$, without loss
of generality we assume that \texttt{prod} is started at epoch $(0,0)$ with
$M^{(0)} = M_1$ and $S_{-1} = 0$.

The last step of epoch $(r,s)$ is the time step $t = t_{(r,s)}$ when either:
\begin{verse}
(C1) $Q^*_t > 4^{S_{r-1}+s} M_t^2$ happens for the first time
\end{verse}
or
\begin{verse}
(C2) $M_t > M^{(r)}$ happens for the first time.
\end{verse}
If epoch $(r,s)$ ends because of (C1),
%YM: why is the next epoch not (r,s+k) where k is the minimal ...
%GS  Sure, this together with a very careful analysis may improve the remainder terms in the bound.
%       However, for simplicity (as we did in Section 3.1 as well!) it's easier to
%       consider that the next epoch is s+1 ... (Otherwise many sums in the proof in the
%       Appendix will be delicate to be written: which values of s within regime r were indeed
%       taken?)
the next epoch is $(r,s+1)$, and the value of $M^{(r)}$
is unchanged. If epoch $(r,s)$ ends because of (C2), the next epoch is $(r+1,0)$,
$S_{r} = S_{r-1} + s$, and $M^{(r+1)} = M_t$.
% We already recalled above that '$M_t$ is defined in Section~\ref{SecUnkM}' (M_t is used in
% the definition of the algorithm). If we write again the formula, then we have to do it above.
% '(Recall that $M_t = \max_{s=1,\ldots,t} \max_{i = 1,\ldots,N} 2^{\lceil \log_2 |x_{i,s}| \rceil}$.)'

Note that within epochs indexed by the same $r$, the payoffs in all
steps but the last one are bounded by $M^{(r)}$. Note also that the
quantities $S_r$ count the number of times an epoch ended because of
(C1). Finally, note that there are $S_{r} - S_{r-1}+1$ epochs
$(r,s)$ for a given $r \geq 0$, indexed by $s = 0, \ldots, S_{r} -
S_{r-1}$.
\begin{theorem}
\label{th:doubling-refined2}
For any sequence of payoffs and for any $n \ge 1$, the cumulative reward of algorithm {\tt prod-MQ} satisfies
\[
\begin{array}{ccccl}
\hX_n & \ge & X_n^* & - & \displaystyle{32 M \sqrt{ q \ln N}} \\
& & & - & \displaystyle{22M \left(1 + \ln N \right) - 2 M \log_2 n - 4 M \lceil (\log_2 \ln N)/2 \rceil} \\
& = & X_n^* & - & \displaystyle{O\left( M \sqrt{q\ln N} + M\ln n + M \ln N\right)}
\end{array}
\]
where
$M = \max_{1\le i \le N}\max_{1 \le t \le n} |x_{i,t}|$ and $\displaystyle{q =
% Had to add this max (1, ...), which does not matter since the quantity
% is still stable by rescaling -- see the last lines of the proof for an explanation
\max \left\{ 1, \, \max_{s \leq n} \frac{Q^*_s}{M_s^2} \right\}}$. \\
\end{theorem}
The proof is in the Appendix.

\section{Second-order bounds for weighted majority}
%GS Important remark:
%      The old bounds of Theorems 4 and 5 were NOT stable by translations,
%      in view of the definition of M as a maximum payoff
%      The solution was however easy: Define M as a maximum range (now denoted by E):
%      E = max_{i,j,t} | x_{j,t} - x_{i,t} | = max_{t=1,\ldots,n} E_t
%      The proofs only needed minor modifications, since we essentially
%      used M to bound variance terms, and can always replace these M by E
%
\label{s:wm}
In this section we derive new regret bounds for the weighted majority
forecaster of \inlinecite{LW94} using a time-varying
learning rate. This allows us to avoid the doubling tricks of
Section~\ref{s:refined} and keep the assumption that no knowledge
on the payoff sequence is available to the forecaster beforehand.

Similarly to the results of Section~\ref{s:refined}, the main term in the
new bounds depends on second-order quantities associated to the sequence of payoffs.
However, the precise definition of these quantities makes the bounds
of this section generally not comparable to the bounds obtained in
Section~\ref{s:refined}.

The weighted majority forecaster using the sequence
$\eta_2,\eta_3,\ldots > 0$ of learning rates assigns at time $t$
a probability distribution $\bp_t$ over the $N$ experts defined
by $\bp_1 = (1/N,\ldots,1/N)$ and
\begin{equation}
\label{EWAdef}
    p_{i,t} = \frac{e^{\eta_t X_{i,t-1}}}{\sum_{j=1}^N e^{\eta_t X_{j,t-1}}}
\qquad
    \mbox{\rm for $i=1,\ldots,N$ and $t \ge 2$.}
\end{equation}
Note that the quantities $\eta_t > 0$ may depend on the past payoffs
$x_{i,s}$, $i = 1,\ldots,N$ and $s = 1,\ldots,t-1$.
The analysis of \inlinecite{AuCeGe02}, for a related variant
of weighted majority, is at the core of the proof of the following lemma
(proof in Appendix).
\begin{lemma}
\label{LmEWA_incr}
Consider
any nonincreasing sequence $\eta_2,\eta_3,\ldots$ of positive learning rates
and any sequence $\bx_1,\bx_2,\ldots \in \R^N$ of payoff vectors.
Define the nonnegative function $\Phi$ by
\begin{eqnarray*}
    \Phi(\bp_t, \, \eta_t, \, \bx_t)
&=&
    - \sum_{i=1}^N p_{i,t} x_{i,t}
    + \frac{1}{\eta_t} \ln {\sum_{i=1}^N p_{i,t} e^{\eta_t x_{i,t}}}
\\ &=&
    \frac{1}{\eta_t} \ln \left(\sum_{i=1}^N p_{i,t} e^{\eta_t (x_{i,t} - \hx_t)} \right)~.
\end{eqnarray*}
Then the weighted majority forecaster~(\ref{EWAdef}) run with the sequence
$\eta_2,$ $\eta_3,\ldots$ satisfies, for any $n \ge 1$ and for any $\eta_1 \geq \eta_2$,
\[
    \hX_n - X_n^*
\ge
    - \left( \frac{2}{\eta_{n+1}} - \frac{1}{\eta_1} \right) \ln N
    - \sum_{t=1}^n \Phi(\bp_t, \, \eta_t, \, \bx_t)~.
\]
\end{lemma}
Let $Z_t$ be the random variable with range $\{x_{1,t},\ldots,x_{N,t}\}$ and distribution $\bp_t$.
Note that $\E Z_t$ is the expected payoff $\hx_t$ of the forecaster using distribution $\bp_t$
at time $t$. Introduce
\[
    \var Z_t
=
\E Z_t^2 - \E^2 Z_t
=
    \sum_{i=1}^N p_{i,t} x_{i,t}^2 - \left( \sum_{i=1}^N p_{i,t} x_{i,t} \right)^2~.
\]
Hence $\var Z_t$ is the variance of the payoffs at time $t$ under the
distribution $\bp_t$ and the cumulative variance $V_n = \var Z_1 + \cdots + \var Z_n$ is the
main second-order quantity used in this section.
The next result bounds $\Phi(\bp_t, \, \eta_t, \, \bx_t)$ in terms of $\var Z_t$.
\begin{lemma}
%GS Reformulated in terms of payoff ranges
\label{LmPhiUB}
For all payoff vectors $\bx_t = (x_{1,t},\ldots,x_{N,t})$,
all probability distributions $\bp_t =  (p_{1,t},\ldots,p_{N,t})$,
and all learning rates $\eta_t \geq 0$, we have
\[
\Phi(\bp_t, \, \eta_t, \, \bx_t) \leq E
%GS We could even have \eta_t E^2/8 right? (By Hoeffding.)
%      Then in the proof, we would have, in view of the value of \eta_t, E^2/(8E_t) \leq E_t/2 instead of E_t...
%      but who cares about a factor 1/2 in the remainder terms? ...
\]
where $E$ is such that
$\abs[x_{i,t} - x_{j,t}] \leq E$ for all $i,j = 1,\ldots,N$.
If, in addition, $0 \le \eta_t \abs[x_{i,t} - x_{j,t}] \le 1$ for all $i,j =1,\ldots,N$,
then
\[
    \Phi(\bp_t, \, \eta_t, \, \bx_t) \le (e-2) \eta_t \var Z_t~.
\]
\end{lemma}
\begin{pf}
The first inequality is straightforward. To prove the second one we use
$e^a \leq 1 + a + (e-2)\,a^2$ for $\abs[a] \le 1$. Consequently,
noting that $\eta_t \abs[x_{i,t} - \hx_t] \le 1$ for all $i$ by assumption, we
have that
\begin{eqnarray*}
\lefteqn{
    \Phi(\bp_t, \, \eta_t, \, \bx_t)
}
\\ &=&
    \frac{1}{\eta_t} \ln \left(\sum_{i=1}^N p_{i,t} e^{\eta_t (x_{i,t} - \hx_t)} \right)
\\ & \le &
    \frac{1}{\eta_t} \ln \left(\sum_{i=1}^N p_{i,t}
    \left( 1 + \eta_t (x_{i,t} - \hx_t) + (e-2) \eta_t^2 (x_{i,t} - \hx_t)^2 \right) \right)~.
\end{eqnarray*}
Using $\ln (1+a) \leq a$ for all $a > -1$ and some simple algebra concludes
the proof of the second inequality.
\qed
\end{pf}
In~\inlinecite{AuCeFrSc02} a very similar result is proven, except that there the variance is further bounded
(up to a multiplicative factor) by the expectation $\hx_t$ of $Z_t$.

%GS Changed things also accordingly in the proof in the rest of the paper.
%      Basically: it suffices to replace all 2M by E! Even in the titles of the subsections...
%
%      Should we say somewhere that what we use, is that the variance is less than E^2/4 (E being
%      a bound on the range)?
\subsection{Known bound on the payoff ranges (E)}
We now introduce a time-varying learning rate based on $V_n$.
For simplicity, we assume in a first time that a bound $E$ on the payoff ranges $E_t$,
defined in (\ref{defER}),
is known beforehand and turn back to the general case in Theorem~\ref{secondordertheoguess}.
The sequence $\eta_2,\eta_3,\ldots$ is defined as
\begin{equation}
\label{goodetaS}
    \eta_t = \min\left\{ \frac{1}{E},\, C \sqrt{\frac{\ln N}{V_{t-1}}} \right\}
\end{equation}
for $t \ge 2$, with $C = \sqrt{2\left(\sqrt{2}-1\right)/(e-2)} \approx
%GS Yishay, I find 1.07 instead of your 0.7 ...
1.07$.

Note that $\eta_t$ depends on the forecaster's past predictions.
This is in the same spirit as the self-confident learning rates considered in~\inlinecite{AuCeGe02}.
\begin{theorem}
\label{secondordertheoguessS}
%GS Payoff ranges here as well
Provided a bound $E$ on the payoff ranges is known beforehand,
i.e., $\max_{t=1,\ldots,n}\max_{i,j=1,\ldots,N} \abs[x_{i,t} - x_{j,t}] \leq E$,
the weighted majority forecaster using the time-varying learning rate~(\ref{goodetaS})
achieves, for all sequences of payoffs and for all $n \geq 1$,
\[
    \hX_n - X_n^*
\ge
%GS Remainder terms changed (decreased) -- see the end of the proof
    - 4 \sqrt{V_n\ln N} - 2 E \ln N - E/2~.
\]
\end{theorem}
\begin{pf}
We start by applying Lemma~\ref{LmEWA_incr} using the learning rate~(\ref{goodetaS}),
and setting $\eta_1 = \eta_2$ for the analysis,
\begin{eqnarray*}
\lefteqn{
    \hX_n - X_n^*
}
\\ & \ge &
    - \left( \frac{2}{\eta_{n+1}} - \frac{1}{\eta_1} \right) \ln N
    - \sum_{t=1}^n \Phi(\bp_t, \, \eta_t, \, \bx_t)
\\* & \ge &
    - 2\max\left\{E\ln N,\, (1/C)\sqrt{V_n\ln N}\right\}
    -(e-2) \sum_{t=1}^n \eta_t \var Z_t
\end{eqnarray*}
where $C$ is defined in~(\ref{goodetaS}) and the second
inequality follows from the second bound of Lemma~\ref{LmPhiUB}.
We now denote by $T$ the first time step $t$ when $V_t > E^2/4$.
Using that $\eta_t \le 1/E$ for all $t$ and $V_T \le E^2/2$, we get
\begin{equation}
\label{thresWM}
    \sum_{t=1}^n \eta_t \var Z_t \le \frac{E}{2} + \sum_{t=T+1}^n \eta_t \var Z_t~.
\end{equation}
We bound the last sum using $\eta_t \leq C \sqrt{(\ln N)/V_{t-1}}$ for $t \geq T+1$
(note that, for $t \geq T+1$, $V_{t-1} \geq V_T > E^2/4 > 0$). This yields
\[
    \sum_{t=T+1}^n \eta_t \var Z_t
\le
    C \sqrt{\ln N} \sum_{t=T+1}^n \frac{V_t - V_{t-1}}{\sqrt{V_{t-1}}}~.
\]
Since $V_t \leq V_{t-1} + E^2/4 $ and $V_{t-1}\geq E^2/4$ for $t\geq T+1$, we have
\begin{eqnarray*}
\frac{V_t - V_{t-1}}{\sqrt{V_{t-1}}} & = & \frac{\sqrt{V_t}+\sqrt{V_{t-1}}}{\sqrt{V_{t-1}}} \, \left( \sqrt{V_t}-\sqrt{V_{t-1}} \right)
\\ & \leq & (\sqrt{2}+1) \left( \sqrt{V_t} - \sqrt{V_{t-1}} \right)
=  \frac{\sqrt{V_t}-\sqrt{V_{t-1}}}{\sqrt{2}-1}~.
\end{eqnarray*}
Therefore, by a telescoping argument,
\begin{eqnarray} \label{varWM}
    \sum_{t=T+1}^n \eta_t \var Z_t
& \le &
    \frac{C\sqrt{\ln N}}{\sqrt{2}-1}\left(\sqrt{V_n} - \sqrt{V_T} \right) \\
\nonumber & \le &
    \frac{C}{\sqrt{2}-1}\sqrt{V_n\ln N}~.
\end{eqnarray}
Putting things together, we have already proved that
\begin{eqnarray*}
\hX_n - X_n^* & \ge & - 2\max\left\{E\ln N,\, (1/C)\sqrt{V_n\ln N}\right\} \\
& & - \frac{e-2}{2} E - \frac{C (e-2)}{\sqrt{2}-1}\sqrt{V_n\ln N}~.
\end{eqnarray*}
In the case when $\sqrt{V_n} \geq C E \sqrt{\ln N}$, the regret $\hX_n - X_n^*$ is bounded from below by
\[
- \left( \frac{2}{C} + \frac{C(e-2)}{\sqrt{2}-1} \right) \sqrt{V_n\ln N} - \frac{e-2}{2} E
\geq - 4\sqrt{V_n\ln N} - E/2~,
\]
where we substituted the value of $C$ and obtained a constant for the leading term equal to
$2\sqrt{2(e-2)}/\sqrt{\sqrt{2}-1} \leq 3.75$.
When $\sqrt{V_n} \le C E \sqrt{\ln N}$, the lower bound is more than
\begin{eqnarray*}
%GS Also changed this line to improve the way we treat the remainder terms
\lefteqn{- 2 E \ln N - \frac{C (e-2)}{\sqrt{2}-1}\sqrt{V_n\ln N} - \frac{e-2}{2} E} \\
& & \ge - 2 E \ln N - 2 \sqrt{V_n \ln N} - E/2~.
\end{eqnarray*}
This concludes the proof.
\qed
\end{pf}

\subsection{Unknown bound on the payoff ranges (E)}
\label{s:E}
We present the adaptation needed when no bound on the real-valued
payoff range is known beforehand. For any sequence of payoff
vectors $\bx_1,$ $\bx_2,\ldots$ and for all $t=1,2,\ldots$,
%GS Estimation of the payoff range now -- %recall from Section~\ref{SecUnkM}
%     The definition of the algorithm, the statement of the theorem,
%     as well as the proof, have been modified.
we define, similarly to Section~\ref{SecUnkM},
a quantity that keeps track of the payoff ranges seen so far.
More precisely, $E_t = 2^k$, where $k \in \mathbb{Z}$ is the smallest
integer such that
$\max_{s=1,\ldots,t} \max_{i,j=1,\ldots,N}$ $\abs[x_{i,s} - x_{j,s}] \le 2^k$.
Now let the sequence $\eta_2,\eta_3,\ldots$ be defined as
\begin{equation}
\label{goodeta}
    \eta_t = \min\left\{ \frac{1}{E_{t-1}},\, C \sqrt{\frac{\ln N}{V_{t-1}}} \right\}
\end{equation}
for $t \ge 2$, with $C = \sqrt{2\left(\sqrt{2}-1\right)/(e-2)}$.

We are now ready to state and prove the main result of this section,
which bounds the regret in terms of the variance of the predictions.
We show in the next section how this bound leads to more intrinsic bounds on the regret.
\begin{theorem}
\label{secondordertheoguess}
Consider the weighted majority forecaster using the time varying learning rate~(\ref{goodeta}).
Then, for all sequences of payoffs and for all $n \geq 1$,
\[
    \hX_n - X_n^*
\ge
%Remainder terms also changed
    - 4 \sqrt{V_n\ln N} - 4 E \ln N - 6 E
\]
where $E = \max_{t=1,\ldots,n}\max_{i,j=1,\ldots,N} \abs[x_{i,t} - x_{j,t}]$.
\end{theorem}
\begin{pf}
The proof is similar to the one of Theorem~\ref{secondordertheoguessS},
we only have to deal with the estimation of the payoff ranges.
We apply again Lemma~\ref{LmEWA_incr},
\begin{eqnarray*}
    \hX_n - X_n^*
& \ge &
    - \left( \frac{2}{\eta_{n+1}} - \frac{1}{\eta_1} \right) \ln N
    - \sum_{t=1}^n \Phi(\bp_t, \, \eta_t, \, \bx_t)
\\ & \ge &
    - 2\max\left\{E_n\ln N,\, (1/C)\sqrt{V_n\ln N}\right\}
    - \sum_{t=1}^n \Phi(\bp_t, \, \eta_t, \, \bx_t)
\\ &=&
    - 2\max\left\{E_n\ln N,\, (1/C)\sqrt{V_n\ln N}\right\}
\\ &&
    - \sum_{t\in\mathcal{T}} \Phi(\bp_t, \, \eta_t, \, \bx_t)
    - \sum_{t\not\in\mathcal{T}} \Phi(\bp_t, \, \eta_t, \, \bx_t)
\end{eqnarray*}
where $C$ is defined in~(\ref{goodeta}),
and $\mathcal{T}$ is the set of time steps $t \geq 2$ when $E_t = E_{t-1}$
(note that $1 \not\in\mathcal{T}$ by definition).
Thus $\mathcal{T}$ is a finite union of intervals of integers,
$\mathcal{T} = \llbracket 1,\,n \rrbracket \setminus \{ t_1,\ldots,t_R \}$,
where we denote $t_1 = 1$ and let $t_2,\ldots,t_R$ be the time rounds $t \geq 2$
such that $E_t \ne E_{t-1}$.

Using the second bound of Lemma~\ref{LmPhiUB} on $t\in\mathcal{T}$
(since, for $t\in\mathcal{T}$, $\eta_t E_t \leq E_t/E_{t-1} = 1$)
and the first bound of Lemma~\ref{LmPhiUB} on $t\not\in\mathcal{T}$,
which in this case reads $\Phi(\bp_t, \, \eta_t, \, \bx_t) \le E_t$,
we get
\begin{eqnarray}
\nonumber
    \hX_n - X_n^*
& \ge &
    - 2\max\left\{ E_n\ln N,\, (1/C)\sqrt{V_n\ln N}\right\}
\\ &&
\label{Bd1}
    -(e-2) \sum_{t \in \mathcal{T}} \eta_t \var Z_t - \sum_{t \not\in \mathcal{T}} E_t~.
\end{eqnarray}
We consider the $r$-th regime, $r=1,\ldots,R$, that is, the time steps $s$
between $t_r +1$ and $t_{r+1}-1$ (with $t_{R+1}=n$ by convention
%GS -- Here also we need to deal with the case when a new regime starts at n
%         Turns out that it is possible to do it in a simpler way than in Section 3.
whenever $t_R < n$).
For all these time steps $s$, $E_s = E_{t_r}$.
We use the same arguments that led to~(\ref{thresWM}) and~(\ref{varWM}):
denote by $T_r$ the first time step $s \geq t_r + 1$ when $V_s > E_{t_r}^2/4$.
Then,
\[
\sum_{s = t_r +1}^{t_{r+1}-1} \eta_t \var Z_t \leq \frac{E_{t_r}}{2} + \frac{C \sqrt{\ln N}}{\sqrt{2}-1} \left( \sqrt{V_{t_{r+1}-1}}
- \sqrt{V_{T_r}} \right)~.
\]
Summing over $r=1,\ldots,R$ and noting that a telescoping argument is given by $V_{t_r} \leq V_{T_r}$,
\[
\sum_{t \in \mathcal{T}} \eta_t \var Z_t \leq \frac{C \sqrt{\ln N}}{\sqrt{2}-1} \sqrt{V_n} + \frac{1}{2} \sum_{r=1}^R E_{t_r}~.
\]
We deal with the last sum (also present in~(\ref{Bd1})) by noting that
\[
    \sum_{t \not\in \mathcal{T}} E_t = \sum_{r=1}^R E_{t_r}
\le
    \sum_{r=-\infty}^{\lceil \log_2 E \rceil} 2^r
\le
    2^{1+\lceil \log_2 E \rceil}
\le
    4 E~.
\]
Putting things together,
\begin{eqnarray*}
    \hX_n - X_n^*
& \ge &
    - 2\max\left\{ E_n\ln N,\, (1/C)\sqrt{V_n\ln N}\right\} \\
& &    -\frac{(e-2)C \sqrt{\ln N}}{\sqrt{2}-1} \sqrt{V_n} - 2 e \, E~.
\end{eqnarray*}
The proof is concluded, as the previous one, by noting that $E_n \leq 2 E$.
\qed
\end{pf}

%NCB Changed title
% \subsection{Randomized prediction and non-expected regret}
\subsection{Randomized prediction and actual regret}
\label{NONexpR}

In this paper, the focus is on improved bounds for the expected regret.
After choosing a probability distribution $\bp_t$ on the actions,
the forecaster gets $\hx_t = x_{1,t}p_{1,t} + \cdots + x_{N,t}p_{N,t}$
as a reward. In case randomized prediction is considered,
after choosing $\bp_t$, the forecaster draws an action $I_t$ at random
according to $\bp_t$ and gets the reward $x_{I_t,t}$, whose
conditional expectation is $\hx_t$. In this version of the
game of prediction, the aim is now to minimize
the (actual) regret, defined as the difference
between $x_{I_1,1} + \cdots + x_{I_n,n}$ and $X_n^*$.

Bernstein's inequality for martingales (see, e.g., \opencite{Fre75})
shows however that the actual regret of any forecaster is bounded by the
expected regret with probability $1-\delta$ up to deviations of the order
of ${\sqrt{V_n \ln (n/\delta)} + M \ln (n/\delta)}$.
These deviations are of the same order of magnitude as the bound of
Theorem~\ref{secondordertheoguess}.
Unless we are able to apply a sharper concentration result than Bernstein's inequality,
no further refinement of the above bounds is worthwhile. In particular,
in view of the deviations from the expectations,
as far as actual regret is concerned,
we may prefer the results of Section~\ref{s:wm} to those
of Section~\ref{s:refined}.
%GS
The next section, as well as Section~\ref{SecAppli}, explain how bounds in terms of $\sqrt{V_n}$
lead to many interesting bounds on the regret that do not depend on
%NCB
% on-line quantities.
quantities related to the forecaster's rewards.

%GS I also moved this section and renamed it.
%      Because actually, it has nothing to do with translations!
%      The quantities $\mu_t$ are indeed chosen in hindsight. And we only use the fact
%      that the variances are less than expected squared values.
%NCB Changed title
%\subsection{Upper bounding the variance $V_n$ of the on-line prediction.}
\subsection{Bounds on the forecaster's cumulative variance}
\label{TranslSec}
In this section we show a first way to deal with the dependency of the bound
%NCB
% in terms of the on-line quantities that arise in
on $V_n$, the forecaster's cumulative variance.
Section~\ref{SecAppli} will illustrate this further.

Recall that $Z_t$ is the random variable which takes the value $x_{i,t}$ with probability $p_{i,t}$,
for $i=1,\ldots,N$. The main term of the bound stated in Theorem~\ref{secondordertheoguess}
contains $V_n = \var Z_1 + \cdots + \var Z_n$. Note that $V_n$ is therefore smaller
than all quantities of the form
\[
    \sum_{t=1}^n \sum_{i=1}^N p_{i,t} \left( x_{i,t} - \mu_t \right)^2
\]
where $(\mu_t)_{t \ge 1}$ is any sequence of real numbers which may be chosen
in \textsl{hindsight}, as it is not required for the definition of the forecaster.
(The minimal value of the expression is obtained for $\mu_t=\hat{x}_t$.)
This gives us a whole family of upper bounds, and we may choose for the analysis the
most convenient sequence of $\mu_t$.

To provide a concrete example,
recall the definition~(\ref{defER}) of payoff effective range $E_t$
and consider the choice $\mu_t = \min_{j=1,\ldots,N} x_{j,t} + E_t/2$.
\begin{corollary} \label{CorRange}
The regret of the weighted majority forecaster with variable learning
rate~(\ref{goodeta}) satisfies
\[
    \hX_n - X_n^*
\ge
    - 2\sqrt{(\ln N)\sum_{t=1}^n E_t^2} -
%GS (see the changes in the previous section)
4 E \ln N - 6 E
\]
where $E$ is a bound on the payoff ranges, $E = \max_{t=1,\ldots,n} E_t$.
\end{corollary}
The bound proposed by Corollary~\ref{CorRange} shows that for an effective range of $E$,
say if the payoffs all fall in $[0,E]$, the regret is lower bounded by a quantity
equal to $- 2 E \sqrt{n \ln N}$ (a closer look at the proof of Theorem~\ref{secondordertheoguess}
shows that this constant factor is less than $1.9$,
% This is because the threshold at $E_{t_r}^2/4$ determining $T_r$ by means of the $Vs$
% has been set quite arbitrarily. A value of $a E^2$, $a \geq 1/4$, would lead to a bound with a smaller
% constant factor in the leading term, at the cost of larger constant terms in the remainder.
and could be made as close to $2 \sqrt{(e-2)}
%GS Added this -- see the discussion about the constants in Section 5.4 why
= \sqrt{2} \sqrt{2 \,(e-2)}$
as desired).
The best leading constant for such bounds is, to our knowledge, $\sqrt{2}$ (see~\opencite{CBL05}).
This shows that the improved dependence in the bound does not come at a significant increase
in the magnitude of the leading coefficient.
%NCB moved these remarks after the statement
When the actual ranges are small, these bounds give a considerable advantage.
Such a situation arises, for instance, in the setting of on-line portfolio selection,
when we use linear upper bound on the regrets (see, e.g., the {\sc eg} strategy by \opencite{HeScSiWa98}).
Moreover, we note that Corollary~\ref{CorRange} improves on a result of \inlinecite{AlNe04},
who show a regret bound, in terms of the cumulative effective range, whose main
term is ${5.7 \sqrt{2 M (\ln N) \sum_{t=1}^n E_t}}$, for a given bound $M$ over the payoffs.

Finally, we note that using translations of payoffs for {\tt prod}-type algorithms, as suggested by
Section~\ref{s:prodtrans}, may be worthwhile as well, see Corollary~\ref{SmOrLg-2} below.
However, unlike the approach presented here for the weighted majority
based forecaster, there the payoffs have to be translated
%GS
explicitly and on-line
by the forecaster, and thus, each translation rule corresponds to a different forecaster.

\subsection{Extension to problems with incomplete information}
\label{s:review}
%
%GS -- New section -- It should be expanded a bit for this paper to be
%         self-contained... (I didn't do it because I first wanted to be sure
%         that you find this section relevant. I definitely do. People from the
%         game-theory community in Paris were very happy to have a statistical
%         explanation of the difference in rates between partial monitoring and bandit.)
%
%YM My feeling is that this sub-section can be omitted. Both its a
%       different model, and the results/remarks require complete knowledge
%       of the previous results. Maybe we can just add a a paragraph to the
%       conclusion section.
%
%GS Alternatively, this section can be shortened. But I think that this
%      is its right place.
%
%NCB I would suggest to keep it. In the worst case the reviewers will ask to drop it.
%
An interesting issue is how the second-order bounds of this
section extend to incomplete information problems. In the
literature of this area, exponentially weighted averages of estimated
cumulative payoffs play a key role (see, for instance, \opencite{AuCeFrSc02}
for the multiarmed bandit problem, \opencite{CeLuSt04} for
label-efficient prediction, and \opencite{PiSc01}, \opencite{CLS04}
for prediction under partial monitoring).

A careful analysis of the proofs therein shows that the order of magnitude of the bound
on the regret is given by the root of the sum of the conditional variances of the
estimates of the payoffs used for prediction,
\[
\sqrt{(\ln N) \sum_{t=1}^n \mathbb{E}_t \left[ \sum_{i=1}^{N}
p_{i,t} \left( \widetilde{x}_{i,t} \right)^2
- \left( \sum_{i=1}^{N}
p_{i,t} \widetilde{x}_{i,t} \right)^2
\right]}~.
\]
Here we denote by $\widetilde{x}_{i,t}$ the (unbiased) estimate available for $x_{i,t}$
(whose form varies depending on the precise setup and the considered strategy),
by $\bp_t = (p_{1,t},\ldots,p_{N,t})$ the probability distributions over the actions,
and by $\mathbb{E}_t$ the conditional expectation with respect to the information available
up to round $t$ (for instance, in multiarmed bandit problems, this information is the past payoffs).
Note that the conditioning in $\mathbb{E}_t$ determines the values of the payoffs
$\bx_t = (x_{1,t},\ldots, x_{N,t})$ and of $\bp_t$.

% As indicated above, the orders of magnitude of the
% main term in an upper bound on the regret, both in $n$ and $N$, are given by
%
In setups with full monitoring, that is, for the setups considered in this paper,
no estimation is needed,
$\widetilde{x}_{i,t} = x_{i,t}$, and the bound is exactly that of Theorem~\ref{secondordertheoguess}.

In multiarmed bandit problems (with payoffs in, say, $[-M,M]$),
the estimators are given by $\widetilde{x}_{i,t} = (x_{i,t}/p_{i,t}) \mathbb{I}_{[I_t = i]}$
where $I_t$ is the index of the chosen component of the payoff vector. Now,
\begin{equation} \label{varest}
\mathbb{E}_t \left[ p_{i,t}\,\widetilde{x}_{i,t}^2 \right] = {x}_{i,t}^2 \leq M^2~.
\end{equation}
Summing over $i=1,\ldots,N$ and $t=1,\ldots,n$ the bound $M \sqrt{n N \ln N}$
of \inlinecite{AuCeFrSc02} is recovered.

In label-efficient prediction problems, $\widetilde{x}_{i,t} = (x_{i,t}/\varepsilon) Z_t$,
where the $Z_t$ are i.i.d.\ random variables distributed according to a Bernoulli
distribution with parameter $\varepsilon \sim m/n$. Then,
\[
\mathbb{E}_t \left[ p_{i,t} \, \widetilde{x}_{i,t}^2 \right] = p_{i,t} \frac{{x}_{i,t}^2}{\varepsilon} \leq
p_{i,t} \frac{M^2}{\varepsilon}~.
\]
Summing over $i=1,\ldots,N$ and $t=1,\ldots,n$ we recover the bound
$M \sqrt{(n/\varepsilon) \ln N} \sim M n \sqrt{(\ln N)/m}$ of \inlinecite{CeLuSt04}.

Finally, in games with partial monitoring, the quantity~(\ref{varest}) is less than
$M^2 t^{-1/3} N^{2/3} (\ln N)^{- 1/3}$. Summing over $i=1,\ldots,N$ and $t=1,\ldots,n$
we recover the $M n^{2/3} N^{2/3} (\ln N)^{1/3}$ bound of \inlinecite{CLS04}.

In conclusion, the faster $\sqrt{n}$ rate in bandit problems, as opposed to the $n^{2/3}$
rate in problems of prediction under partial monitoring, is due to better statistical
performances (i.e., smaller conditional variance) of the available estimators.

%GS Since all this section is now about translations...
\section{Using translations of the payoffs}
\label{SecAppli}
We now consider the bounds derived from those of
Sections~\ref{s:refined} and~\ref{s:wm} in the case when
translations are performed on the payoffs
(Section~\ref{s:prodtrans}). We show that they lead to several
improvements or extensions of earlier results
(Section~\ref{s:SorL}) and also
%YM-2-2-06
%allow the forecaster to avoid any preliminary manipulation on the payoffs (Section~\ref{s:stab}).
relieve the forecaster from the need of any preliminary
manipulation on the payoffs (Section~\ref{s:stab}).

\subsection{On-line translations of the payoffs}
\label{s:prodtrans}
%YM I just moved the sub-section here from Section 3. It has a
%       significant overlap with the next sub-section.
%GS  OK, I modified the introduction (next two paragraphs) of this section.
%       No overlap any more.
%
Note that any on-line forecasting strategy may be used by a meta-forecaster
which, before applying the given strategy, may first translate the payoffs
according to a prescribed rule that may depend on the past.
More formally, the meta-forecaster runs the strategy with the payoffs
$r_{k,t} = x_{k,t} - \mu_t$, where $\mu_t$ is any quantity
possibly based on the past payoffs $x_{i,s}$, for $i = 1,\ldots,N$
and $s = 1,\ldots,t$.

The forecasting strategies of Section~\ref{s:wm} (and the obtained bounds)
are invariant by such translations. This is however not the case for
the {\tt prod}-type algorithms of Section~\ref{s:refined}.
An interesting application is obtained in Section~\ref{s:SorL}
by considering $\mu_t = \hx_t$ where we recall that $\hx_t = x_{1,t}
p_{1,t} + \cdots + x_{N,t} p_{N,t}$ is the forecaster's reward at
time $t$. As the sums $\mu_1 + \cdots + \mu_n$ cancel out in the
difference $\hX_n - X_{k,n}$, we obtain the following corollary
of Theorem~\ref{th:doubling-refined}.
%GS To get a bound stable by translations, it has to be expressed
%      not in terms of an upper bound on the payoffs, but of an upper bound on the payoff ranges.
Note that the remainder term here is now expressed in terms of the effective
ranges~(\ref{defER}) of the payoffs.
\begin{corollary}
\label{cor:doubling-refined2} Given $E > 0$, for all $n \geq 1$ and
all sequences of payoffs
with effective ranges $E_t$ bounded by $E$, the cumulative reward of
algorithm {\tt prod-Q}$(E)$ run using translated payoffs $x_{k,t} -
\hx_t$ satisfies
\begin{eqnarray*}
    \hX_n
\ge
    X_n^* & - & 8 \sqrt{ (\ln N) \max_{s \leq n} R^*_s} \\*
    & - & 2\, E \Bigl( 1 + \log_4 n + 2 \bigl( 1 + \lfloor (\log_2 \ln N)/2 \rfloor \bigr) \ln N \Bigr)~.
\end{eqnarray*}
where the $R^*_s$ are defined as follows. For $1 \leq t \leq n$ and
$k = 1,\ldots,N$, $R_{k,t} = (x_{k,1} - \hx_1)^2 + \cdots + (x_{k,t}
- \hx_t)^2$ and $R^*_t = R_{k^*_t,t}$, where $k^*_t$ is the index of
the action achieving the best cumulative payoff at round $t$ (ties
are broken by choosing the action $k$ with smallest associated
$R_{k,t}$).
\end{corollary}
\begin{remark} \label{OneSideTR}
In one-sided games, for instance in gain games, the forecaster
has always an incentive to translate the payoffs by the minimal
payoff $\mu_t$ obtained at each round $t$,
\[
    \mu_t = \min_{j=1,\ldots,N} x_{k,t}~.
\]
%YM-2-2-06
%This just because,
This is since
for all $j$ and $t$, $(x_{j,t} - \mu_t)^2 \leq
x_{j,t}^2$ in a gain game. The
%YM-2-2-06
%matter
issue
is not so clear however for signed games, and it may be a
delicate issue to determine beforehand if the payoffs should be
translated, and if so, which translation rule should be used. See
also Section~\ref{TranslSec}, as well as Section~\ref{s:SorL}.
\end{remark}

%NCB Changed title
% \subsection{Improvements for small or large payoffs in one-sided (and signed) games.}
\subsection{Improvements for small or large payoffs}
\label{s:SorL}
%GS I changed the introduction of this section as well (to motivate why
%      we need improvements for small or large payoffs).
%
As recalled in Section~\ref{s:firstbd}, when all payoffs have the same sign \inlinecite{FS97}
first showed that Littlestone and Warmuth's weighted majority algorithm~\shortcite{LW94}
can be used to construct a forecasting strategy achieving a regret of order
$\sqrt{M |X_n^*|\ln N} + M\ln N$, where $N$ is the number of actions,
$M$ is a known upper bound on the magnitude of payoffs
($|x_{i,t}| \le M$ for all $t$ and $i$), and $\abs[X_n^*]$ is
the absolute value of the cumulative payoff of the best action (i.e., the
largest cumulative payoff in a gain game or the smallest cumulative loss in a loss game),
see also \inlinecite{AuCeGe02}.

This bound is good when $\abs[X_n^*]$ is small in the one-sided game;
that is, when the best action has a small gain (in a gain game)
or a small loss (in a loss game).
However, one often expects the best expert to be
effective (for instance, because we have many experts and at least one of them is accurate).
%YM-2-2-06
%A
An
effective expert in a loss game suffers a small cumulative loss,
but in a gain game, such an expert should get a large cumulative
payoff $X_n^*$. To obtain a bound that is good when $\abs[X_n^*]$ is
large one could apply the translation $x'_{i,t} = x_{i,t} - M$ (from
gains to losses) or the translation $x'_{i,t} = x_{i,t} + M$ (from
losses to gains). In both cases one would obtain a bound of the form
$\sqrt{M(Mn - |X_n^*|)\ln N}$, which is now suited for effective
experts in gain games and poor experts in loss games, but not for
effective experts in loss games and poor experts in gain games.
Since the original bound is not stable under the operation of
conversion from one type of one-sided game into the other, the
forecaster has to guess whether to play the original game or its
translated version, depending on his beliefs on the quality of the
experts and on the nature
%YM-2-2-06
%(losses or gains) of the game.
of the game (losses or gains).

In Corollary~\ref{SmOrLg-2} we use the sharper bound of Corollary~\ref{cor:doubling-refined2}
to prove a (first-order) bound of the form
\[
\sqrt{M\min\bigl\{|X_n^*|,\, Mn - |X_n^*| \bigr\}\ln N}~.
\]
This is indeed an improvement for small losses or large gains, though
it requires knowledge of $M$. However,
%YM-2-2-06: forward pointer
%Remark~\ref{SmOrLg-withoutM} indicates
in Remark~\ref{SmOrLg-withoutM} we will indicate
how to extend this result to
the case when $M$ is not known beforehand. Note that the
(second-order) bound of Corollary~\ref{SmOrLg} also yields the same
result without any preliminary knowledge of $M$.

We thus recover an earlier result by \inlinecite{AlNe04}. They
proved, in a gain game, for a related
algorithm, and with the previous knowledge of a bound $M$ on the payoffs,
a bound whose main term is
$11.4 \sqrt{M} \min\!\left\{\sqrt{X_n^*},\, \sqrt{M n - X_n^*}\right\}$.
That algorithm was specifically designed to ensure a regret bound of this form,
and is different from the algorithm whose performance we discussed before the
statement of Corollary~\ref{CorRange},
whereas we obtain
the improvements for small losses or large gains as
corollaries of much more general bounds that have other consequences.

\subsubsection{Analysis for exponentially weighted forecasters}
The main drawback of $V_n$, used in Theorem~\ref{secondordertheoguess}, is that it is defined
directly in terms of the forecaster's distributions $\bp_t$.
We now show how this dependence could be removed.
\begin{corollary}
\label{SmOrLg}
Consider the weighted majority forecaster run with the time-varying learning rate~(\ref{goodeta}).
Then, for all sequences of
payoffs in a one-sided game (i.e., payoffs are all nonpositive or all nonnegative),
\[
    \hX_n
\ge
   X_n^*
  - 4\sqrt{|X^*_n| \left( M - \frac{|X_n^*|}{n} \right) \ln N} \\
%GS 75 updated to 39 in view of the improvements in Section 4
  - 39 \,M \max \left\{ 1, \ln N \right\}
\]
where
%GS The notions of payoff range and bound on the payoffs are the same now (one-sided games)
$M = \max_{t=1,\ldots,n} \max_{i=1,\ldots,N} |x_{i,t}|$.
\end{corollary}
\begin{pf}
We give the proof for a gain game. Since the payoffs are in $[0,M]$, we can write
\begin{eqnarray*}
    V_n
& \le &
    \sum_{t=1}^n \left( M \sum_{i=1}^N p_{i,t} x_{i,t}  - \left( \sum_{i=1}^N p_{i,t} x_{i,t} \right)^2 \right)
    = \sum_{t=1}^n ( M - \hx_{t})\hx_t
\\ & \le &
    n \left( \frac{M \hX_n}{n}  - \left( \frac{\hX_n}{n} \right)^2 \right)
=
    \hX_n \left( M - \frac{\hX_n}{n} \right)
\end{eqnarray*}
where we used the concavity of $x \mapsto M x - x^2$.
Assume that $\hX_n \le X_n^*$ (otherwise the result is trivial).
Then, Theorem~\ref{secondordertheoguess} ensures that
\[
    \hX_n - X_n^*
\ge
    -  4\sqrt{X^*_n \left( M - \frac{\hX_n}{n} \right) \ln N} - \kappa
\]
where
%GS Updated \kappa, recomputed the remainder term in the bound of the theorem
$\kappa = 4 M\ln N + 6 M$.
% We indeed solve for $\hX_n$, by introducing $\hL_n = Mn - \hX_n$ and $L_n^* = M n - X_n^*$
We solve for $\hX_n$ obtaining
\[
    \hX_n - X_n^*
\ge
     - 4\sqrt{X^*_n \left( M - \frac{X_n^*}{n}  + \frac{\kappa}{n} \right) \ln N} - \kappa
- 16 \frac{X_n^*}{n} \ln N~.
\]
Using the crude upper bound $X^*_n/n \le M$ and performing some simple algebra,
we get the desired result.
\qed
\end{pf}
%
%      Here is the summary:
%      Optimal bound (needs the knowledge of n) : \sqrt{n/2}
%      Additional factor of 2 to deal with n unknown
%      Additional factor of 2 if n is to be replaced by X_n^* or L_n^*
%      Additional (small!) factor of \sqrt{2 (e-2)} \sim 1.19 if we ignore the sign of the payoffs
%
Similarly to the remark about constant factors in
%YM-2-2-06
% the previous section,
Section~\ref{TranslSec}
the factor $4$ in Corollary~\ref{SmOrLg} can be made as close as
desired to $4 \sqrt{e-2}
%GS I wrote something similar in Section 5.3, the
%      previous %GS explains why I write constants this way
= 2\sqrt{2} \, \sqrt{2 \, (e-2)}$,
which is not much larger than the best known leading constant for improvements for small
losses, $2 \sqrt{2}$, see \inlinecite{AuCeGe02}. But here, we have in addition
an improvement for large losses, and deal with unknown ranges $M$.
(Note, similarly to the discussion in Section~\ref{TranslSec}, the presence of the
same small factor $\sqrt{2 \, (e-2)} \approx 1.2$.)

\subsubsection{Analysis for {\tt prod}-type forecasters}
Quite surprisingly, a bound of the same form as the one shown in
Corollary~\ref{SmOrLg} can be derived from Corollary~\ref{cor:doubling-refined2}.
\begin{corollary}
\label{SmOrLg-2}
Given $M > 0$,
for all $n \geq 1$ and all sequences of payoffs bounded by $M$, i.e.,
$\max_{1 \leq i \leq N} \max_{1 \leq t \leq n} \abs[x_{i,t}] \leq M$,
the cumulative reward of algorithm
%GS {\tt prod-Q}$(M)$ simply replaced by
{\tt prod-Q}$(2M)$,
%      The proof doesn't have to be modified, this was but a typo
run using translated payoffs $x_{k,t} - \hx_t$
in a one-sided game, is larger than
\begin{eqnarray*}
    \hX_n
\ge
   X_n^*
   & - & 8 \sqrt{ 2 M \min \left\{ X_n^*, \, M n - X_n^* \right\} \ln N } \\
& - & 128 \,M \ln N - \kappa - 8\sqrt{2 M (\ln N) \kappa}
\end{eqnarray*}
where
\begin{eqnarray*}
\kappa & = & 4\, M \Bigl( 1 + \log_4 n + 2 \bigl( 1 + \lfloor (\log_2 \ln N)/2 \rfloor \bigr) \ln N \Bigr) \\
& = & \Theta\Bigl( M (\ln n)+M (\ln N) (\ln \ln N) \Bigr)~.
\end{eqnarray*}
\end{corollary}
\begin{pf}
As in the proof of Corollary~\ref{SmOrLg}, it suffices to give the proof for a gain game.
In fact, we apply below the bound of Corollary~\ref{cor:doubling-refined2}, which
is invariant under the change $\ell_{i,t} = M - x_{i,t}$
that converts bounded losses into bounded nonnegative payoffs.

The main term in the bound of Corollary~\ref{cor:doubling-refined2}, with the notations therein,
involves
\begin{equation} \label{EqThref2}
\max_{s \leq n} R^*_s \leq \min \left\{ M \left( X_n^* + \hX_n \right), \,
M \left( 2\,M n - X_n^* - \hX_n \right) \right\}~.
\end{equation}
Indeed, using that $(a-b)^2 \leq a^2 + b^2$ for $a,\,b \geq 0$, we get on the one hand,
for all $1 \leq s \leq n$,
$$R_s^* \leq \sum_{t=1}^s x_{k^*_s,t}^2 + \hx_s^2 \leq M \left( X_{k^*_s,s} + \hX_s \right)
\leq M \left( X_n^* + \hX_n \right)$$
whereas on the other hand, the same techniques yield
\begin{eqnarray*}
R_s^* & = & \sum_{t=1}^s \left( \Bigl(M - x_{k^*_s,t}\Bigr) -  \Bigl(M - \hx_s^2\Bigr) \right)^2 \\*
& \leq & M \left( \Bigl( M s - X_s^* \Bigr) + \Bigl( M s - \hX_s \Bigr) \right)~.
\end{eqnarray*}
Now, we note that for all $s$, $X^*_{s+1} \leq X^*_s + M$, and similarly,
$\hX_{s+1} \leq \hX_{s} + M$. Thus we also have $\max_{s \leq n} R_s^* \leq
M \, (2Mn - X_n^* - \hX_n)$.

Corollary \ref{cor:doubling-refined2}, combined with~(\ref{EqThref2}), yields
\[
\hX_n \geq \hX_n^* - 8 \sqrt{ M (\ln N) \min \left\{ \left( X_n^* + \hX_n \right), \,
\left( 2\,M n - X_n^* - \hX_n \right) \right\} } - \kappa
\]
where $\kappa = 4\, M \Bigl( 1 + \log_4 n + 2 \bigl( 1 + \lfloor (\log_2 \ln N)/2 \rfloor \bigr) \ln N \Bigr)$.
Without loss of generality, we may assume that $\hX_n \leq X_n^*$ and get
\[
\hX_n \geq \hX_n^* - 8 \sqrt{ 2 M (\ln N) \min \left\{ X_n^*, \,
\left( M n - \hX_n \right) \right\} } - \kappa~.
\]
Solving for $\hX_n$ and performing simple algebra in case the minimum
is achieved by the term containing $\hX_n$ concludes the proof.
\qed
\end{pf}
%
%GS
\begin{remark}
\label{SmOrLg-withoutM}
The forecasting strategy of Theorem~\ref{th:doubling-refined2},
when used by a meta-forecaster translating the payoffs by $\hx_t$,
achieves an improvement for small or large payoffs of the form
\[
M \sqrt{\min \left\{ \max_{s \leq n} \frac{X_s^*}{M_s}, \, \max_{s \leq n} \frac{s M_s - X_s^*}{M_s} \right\} }
\]
without previous knowledge of $M$.
\end{remark}

\subsubsection{The case of signed games}
\label{r:smorlg-signed}
%YM Need to redo this remark, given the changes in section 2
%GS Corrected.
%
The proofs of Corollaries~\ref{SmOrLg} and~\ref{SmOrLg-2} reveal
that the assumption of one-sidedness cannot be relaxed. However, we
may also prove a version of the improvement
for small losses or for large gains suited to signed games. Remember
that, as explained in Section~\ref{s:firstsigned}, a meta-forecaster
may always
convert a signed game into a one-sided game by performing
a suitable translation on the payoffs, and then apply a strategy
for one-sided games. Since Corollary~\ref{cor:doubling-refined2} and
Theorem~\ref{secondordertheoguess} are stable under general
translations, applying them to the payoffs $x_{i,t}$ or to a
translated version of them $x'_{i,t}$ results in the same bounds. If
the translated version $x'_{i,t}$ correspond to a one-sided game,
then the bounds of Corollaries~\ref{SmOrLg} and~\ref{SmOrLg-2} may
be applied. Using $x'_{i,t} = x_{i,t} - \min_{j=1,\ldots,N} x_{j,t}
\geq 0$ and $x'_{i,t} = x_{i,t} - \max_{j=1,\ldots,N} x_{j,t} \leq
0$ for the analysis, we may show, for instance, that for any signed
game the forecaster of Theorem~\ref{secondordertheoguess}
ensures that the regret is bounded by a quantity whose main term is less
than
\begin{eqnarray*}
\min & \displaystyle{\left\{ \sqrt{(\ln N) \max_{j=1,\ldots,N}
\left( \sum_{t=1}^n \left( x_{j,t} - \min_{i=1,\ldots,N} x_{i,t} \right) \right)}~, \right.} \\*
& \displaystyle{\left. \sqrt{(\ln N) \min_{j=1,\ldots,N}
\left( \sum_{t=1}^n \left( \max_{i=1,\ldots,N} x_{i,t} - x_{j,t} \right) \right)} \ \right\}~.}
\end{eqnarray*}
%YM The above bound seems incorrect. Is it missing only a square?!
%GS  No, it's correct (but I agree it could be stated in a stronger version,
%       I mostly chose a readable version).
This bound is obtained without any previous knowledge of a bound $M$
on the payoffs, and is sharper than both bounds~(\ref{EqTr1Sec2})
and~(\ref{EqTr2Sec2}).
It may be interpreted as an improvement for small or large cumulative payoffs.

%GS Tried to summarize quickly the discussion we had about stability in the
%      older versions of this paper. Feel free to shorten it if you think it's needed.
%      But since this is now at the very end of the paper, I think it can be a bit philosophical
%      and formulate some general remarks rather than give a neat conclusion.
%
%      Most importantly, this section explains the form of most of our bounds, e.g.,
%      why the remainder terms of some bounds involve the effective range
%      $E$ rather than the bound $M$ on the magnitude of the payoffs.
%
\subsection{What is a ``fundamental'' bound?}
\label{s:stab}
%GS This section is not in final form. Please help me getting it clear.
%NCB I did some rewriting
%GS To me, it's OK now.
Most of the known regret bounds are not stable under natural transformations
%NCB Put comment as footnote.
of the payoffs, such as translations and rescalings.\footnote{
Here we do not distinguish between stable bounds and stable algorithms
because all the stability properties we consider for the bounds are due
to a corresponding stability of the prediction scheme they are derived from.
When a stable algorithm does not achieve a stable bound, it suffices to
optimize the bound in hindsight, thanks to the stability properties of the
prediction scheme.
}
If a regret bound is not stable, then a
(meta-)prediction  algorithm might be willing
to manipulate the payoffs in order to achieve a better regret.
However, in general it is hard to choose the payoff
transformation that is best for a given and unknown sequence of
payoffs.  For this reason, we argue that regret bounds that are stable
under payoff transformations are, in some sense, more fundamental than others.
The bounds that we have derived in this paper are based on sums of squared payoffs.
They are not only generally tighter than the previously known bounds, but also stable
under different transformations, such as those described below
(in what follows, we use $x_{i,t}'$ to indicate a transformed payoff).
\subsubsection*{Additive translations: $\quad x'_{i,t} = x_{i,t} - \mu_t$.} Note
that the regret (of playing a fixed sequence $\bp_1,\,\bp_2,\ldots$)
is not affected by this transformation. Hence, stable bounds should
not change when payoffs are translated. As already explained in
Section~\ref{s:SorL}, translations can be used to turn a gain game
into a loss game and vice versa.
%NCB deleted
% We recall that we may turn a gain game with payoff range $[0,M]$
% into a loss game with payoff range $[-M,0]$ by performing the
% translation $x'_{i,t} = x_{i,t} - M$.
% (This particular case of translation will be referred to as \emph{conversion}.
% The bounds of Section~\ref{s:SorL} are stable by conversion.)
%
% More generally, any signed game can be transformed in an
% equivalent gain or loss game without knowing the payoff ranges,
% as explained in Section~\ref{r:smorlg-signed}.
% For instance, to transform the payoff $x_{i,t}$ of a signed game
% in the payoff $x'_{i,t}$ of a gain game we just compute
% $x'_{i,t} = x_{i,t}-\min_{j=1,\ldots,N} x_{j,t}$.

The invariance by general translations is the hardest to obtain,
and this paper is the first one to show tight translation-invariant
bounds that depend on the specific sequence of payoffs rather than just
on its length (see Corollary~\ref{cor:doubling-refined2},
Theorem~\ref{secondordertheoguess} and some of their corollaries,
e.g., Corollary~\ref{CorRange}).
It is also important to remark that, in a stable bound, not only the leading term,
but also the smaller order terms, have to be stable under translations.
This is why the smaller order terms of Corollary~\ref{cor:doubling-refined2}
and Theorem~\ref{secondordertheoguess} involve bounds on the payoff ranges
$x_{i,t} - x_{j,t}$ rather than just on the payoffs $x_{i,t}$.

\subsubsection*{Rescalings: $\quad x'_{i,t} = \alpha\,x_{i,t}$, $\alpha > 0$.}
As this transformation causes the regret to be multiplied by a factor of $\alpha$,
stable bounds should only change by the same factor $\alpha$.
Obtaining bounds that are stable under rescalings is not always easy
when the payoff ranges are not known beforehand, or when we try to
get bounds sharper than the basic zero-order bounds discussed in
Section~\ref{s:zerobd}.
%
%GS ... Give one or two examples of undesirable terms? The difficulties are
%      often in the remainder terms (except in the bounds of a previous version of this paper,
%      were even some main terms were not stable!). We could write:
%
For instance, the application of a doubling trick on the magnitude of the payoffs,
or even the use of more sophisticated incremental techniques,
may lead to small but undesirable $M \ln (Mn)$ terms, which behave badly upon rescalings.
This was the case with the remainder term $M \ln (1+\abs[X_n^*])$ in Theorem~2.1
by \inlinecite{AuCeGe02} where they assume knowledge of the payoff range but seek sharper bounds.
%NCB Deleted
% Other examples of bounds that are not stable under rescalings
% are given in the main theorems of
% a previous version of this paper, \inlinecite{CeMaSt05}, which involve
% main terms of the form $\sqrt{ \max \{ 1, \, \, B_n \} \ln N}$, where $B_n$ are
% terms that are scalable, that is, they are multiplied by $\alpha^2$
% when the payoffs are multiplied by $\alpha$.
%

Note also that forecasters with scaling-invariant bounds should require no previous
knowledge on the payoff sequence (such as the payoff range) as this information is
scale-sensitive. This is why, for instance, the bounds of Theorems~\ref{th:doubling-refined}
and~\ref{secondordertheoguessS} cannot be considered scaling-invariant.
However, modifications of these forecasters that increase their adaptiveness lead to
Theorems~\ref{th:doubling-refined2} and~\ref{secondordertheoguess}. There we could derive
scaling-invariant bounds by using forecasters based on updates which are defined in terms
of quantities that already have this type of invariance.

\medskip

%GS Maybe too personal?
%NCB I tried to rewrite it
% Our feeling is that all bounds in the literature should be stable by
% rescaling, this is only a matter of an accurate enough analysis. The
% stability by translations was however a greater challenge since we wanted to obtain
% sharpness at the same time (otherwise, the zero-order bounds would have done the job).
%
Whereas translation-invariant bounds that are also sharp are generally hard to obtain,
we feel that any bound can be made stable with respect to rescalings via a reasonably
accurate analysis.

Unstable bounds can lead the meta-forecaster to Cornelian dilemmas.
Consider for the instance the bound~(\ref{EqTr3Sec2}) by \inlinecite{AlNe04}.
If we use a meta-forecaster that translates payoffs by a quantity $\mu_t$
(possibly depending on past observations), then the bound takes the form
\[
\sqrt{M (\ln N) \max_{t=1,\ldots,n} \sum_{s=1}^t \abs[x_{k_t^*,s} - \mu_s]} + M\ln N~.
\]
Note that the choice $\mu_t = -M$ (or $\mu_t = \min_{j=1,\ldots,N} x_{j,t}$)
yields the improvement for small payoffs~(\ref{EqTr1Sec2}) and the choice
$\mu_t = M$ (or $\mu_t = \max_{j=1,\ldots,N} x_{j,t}$) yields
the improvement for large payoffs~(\ref{EqTr2Sec2}).
In general, the above bound is tight if, for a large number of rounds,
all payoffs $x_{j,t}$ at a given round $t$ are close to a common value,
and we may guess this value to choose $\mu_t$ accordingly.
In Section~\ref{r:smorlg-signed}, on the other hand, we show that
Corollaries~\ref{SmOrLg} and~\ref{SmOrLg-2} propose bounds that need no
preliminary choices of $\mu_t$ and are better than both~(\ref{EqTr1Sec2}) and~(\ref{EqTr2Sec2}).

\section{Discussion and open problems}
%GS This section probably needs some work.
%NCB I have done some rewriting.
%GS To me, this also is OK.
We have analyzed forecasting algorithms that work indifferently in loss games,
gain games, and signed games. In Corollary~\ref{cor:doubling-refined2} and
Theorem~\ref{secondordertheoguess} we have shown, for these forecasters, sharp
regret bounds that are stable under rescalings and general translations.
These bounds lead to improvements for small or large payoffs in one-sided games
(Corollaries~\ref{SmOrLg} and~\ref{SmOrLg-2}) and do not assume any preliminary
information about the payoff sequence.\footnote{
Whereas the bound of Theorem~\ref{secondordertheoguess} is already stated this way,
we recall that it is easy to modify the forecaster used to prove
Corollary~\ref{cor:doubling-refined2} in order to dispense with the need of any
preliminary knowledge of a bound $E$ on the payoff ranges.}

A practical advantage of the weighted majority forecaster is that
its update rule is completely incremental and never needs to reset
the weights. This in contrast to the forecaster {\tt prod-MQ} of
Theorem~\ref{th:doubling-refined2} that uses a nested doubling
trick. On the other hand, the bound proposed in
Theorem~\ref{secondordertheoguess} is not in closed form, as it
still explicitly depends through $V_n$ on the forecaster's rewards
$\hx_t$. We therefore need to solve for the regrets as we did, for
instance, in Sections~\ref{TranslSec} and~\ref{s:SorL}. Finally, it
was also noted in Section~\ref{TranslSec} that the weighted majority
forecaster update is invariant under translations of the payoffs.
This is not the case for the {\tt prod}-type forecasters, which need
to perform translations explicitly. Though in general it may be
difficult to determine beforehand what a good translation could be,
Corollaries~\ref{cor:doubling-refined2} and~\ref{SmOrLg-2}, as well
as Remark~\ref{OneSideTR}, indicate some general effective
translation rules.

Several issues are left open:
\begin{itemize}
\item[--] Design and analyze incremental updates for the {\tt prod}-type forecasters of Section~\ref{s:refined}.
\item[--] Obtain second-order bounds with updates that are not multiplicative; for instance, updates based on
the polynomial potentials (see \opencite{CeLu03}).
These updates could be used as basic ingredients to derive forecasters achieving optimal orders of magnitude
on the regret when applied to problems such as nonstochastic multiarmed bandits, label-efficient prediction,
and partial monitoring.
%GS Updated in view of the new Section 5.2
Note that, to the best of our knowledge, in the literature about incomplete information problems only
exponentially weighted averages have been able to achieve these optimal rates
(see Section~\ref{s:review} and the references therein).
%
%for instance, \opencite{AuCeFrSc02}, \opencite{CeLuSt04}, \opencite{PiSc01}, \opencite{CLS04}).
\item[--] Extend the analysis of {\tt prod}-type algorithms to obtain an oracle inequality of the form
\[
    \hX_n \ge \max_{k=1,\ldots,N} \left( X_{k,n} - \gamma_1 \sqrt{Q_{k,n} \ln N} \, \right) - \gamma_2 M \ln N
\]
where $\gamma_1$ and $\gamma_2$ are absolute constants. Inequalities of this form can be viewed as game-theoretic
versions of the model selection bounds in statistical learning theory.
\end{itemize}

\appendix

\section*{Proof of Theorem~\ref{th:doubling-refined2}}

We use some additional notation for the proof:
$(r,s)-1$ denotes the epoch right before $(r,s)$;
that is, $(r,s-1)$ when $s > 0$, and $(r-1,S_{r-1} - S_{r-2})$
when $s = 0$.
For notational convenience,
$t_{(0,0)-1}$ is conventionally set to $0$.

\begin{pf}
The proof combines the techniques from
Theorems~\ref{th:doubling-refined} and~\ref{guess-M}.
As in the proof of Theorem~\ref{guess-M}, we denote by $(R,S_R-S_{R-1})$
the index of the last epoch and let $t_{(R,S_R-S_{R-1})} = n$.

We assume $R \ge 1$ and $S_R \geq 1$.
Otherwise, if $R = 0$, this means that $M_t = M^{(0)}$
for all $t \leq n-1$, and the strategy, and thus the proposed bound,
reduces to the one of Theorem~\ref{th:doubling-refined}.
The case $S_R = 0$ is dealt with at the end of the proof.
In particular, $S_R \geq 1$ implies that some epoch ended at time $t$ when $Q^*_t > 4^{S_R-1} M_t^2$.
This implies that $q \ge 4^{S_R-1} (\geq 1)$, which in turn implies
$2^{S_R} \leq 2 \sqrt{q}$ and $S_R \le 1 + (\log_2 q)/2$.

%YM: do we need $M^{(R+1)}$, maybe we can just stay with $M_n$.
%GS We need it in many summations below (among others: in the definition of \kappa)
%      Introducing $M^{(R+1)}$ is by far more convenient
Denote $M^{(R+1)} = M_n$.
Note that at time $n$ we have either $M_n \le M^{(R)}$, implying $M_n = M^{(R+1)} = M^{(R)}$,
or we have $M_n > M^{(R)}$, implying $M_n = M^{(R+1)} = 2 M^{(R)}$.
In both cases, $M^{(R)} \le M^{(R+1)} \le 2M$.
Furthermore, $M^{(s)} \ge 2^{s-r} M^{(r)}$ for each $0 \leq r \leq s \leq R$,
and thus~(\ref{sum2power}) holds for $s \leq R$ with $M_{t_r}$ replaced by $M^{(r)}$.

Similar to the proof of Theorem~\ref{th:doubling-refined}, for each epoch $(r,s)$, let
\begin{eqnarray*}
    X_k^{(r,s)} = \!\!\! \sum_{t=t_{{(r,s)}-1}+1}^{t_{(r,s)}-1} \!\!\! x_{k,t}~,
\quad
    Q_k^{({r,s})} = \!\!\! \sum_{t=t_{{(r,s)}-1}+1}^{t_{(r,s)}-1} \!\!\! x_{k,t}^2~,
\quad
    \hX^{({r,s})} = \!\!\! \sum_{t=t_{{(r,s)}-1}+1}^{t_{(r,s)}-1} \!\!\! \hx_{t}
\end{eqnarray*}
where the sums are over all the time steps $t$ in epoch $(r,s)$ except the last one, $t_{(r,s)}$.
We also denote $k_{(r,s)} = k^*_{t_{(r,s)} - 1}$ the index of the best overall expert up to time $t_{(r,s)}-1$
(one time step before the end of epoch $(r,s)$).

We upper bound the cumulative payoff of the best action as
\begin{equation}
\label{Eqinduc}
X_n^* \leq \sum_{r=0}^R \left( M^{(r+1)} + (S_r - S_{r-1}) M^{(r)} +
\sum_{s=0}^{S_{r} - S_{r-1}} X_{k_{(r,s)}}^{(r,s)} \right)
\end{equation}
by using the same argument by induction as in~(\ref{EqCumDT2}).
More precisely, we write, for each $(s,r)$,
\begin{eqnarray*}
X_{k_{(r,s)}, t_{(r,s)}-1}
& = &
X_{k_{(r,s)}}^{(r,s)} + M_{t_{(r,s)-1}} + X_{k_{(r,s)-1}, t_{(r,s)-1}-1} \\
& \leq &
X_{k_{(r,s)}}^{(r,s)} + M_{t_{(r,s)-1}} + X_{k_{(r,s)-1}, t_{(r,s)-1}-1}~.
\end{eqnarray*}
We note that $M_{t_{(r,s)-1}} = M^{(r)}$ whenever $0 \leq s < S_r - S_{r-1}$ and
$M_{t_{(r,s)-1}} = M^{(r+1)}$ otherwise.
This and
\[
X_n^* \leq X_{n-1}^* + M^{(R+1)} = X_{k_{(R,S_R)}, t_{(R,S_R)}-1} + M^{(R+1)}
\]
show~(\ref{Eqinduc}) by induction.

Let
\[
\kappa = \sum_{r=0}^R \left( M^{(r+1)} + (S_r - S_{r-1}) M^{(r)} \right)~.
\]
To show a bound on $\kappa$ note that~(\ref{sum2power}) implies
\begin{equation}
\label{sumR}
    \sum_{r=0}^R M^{(r+1)} \le 2 M^{(R)} + M^{(R+1)} \le 3 M^{(R+1)} \leq 6M
\end{equation}
and
\[
\sum_{r=0}^R (S_r - S_{r-1}) M^{(r)} \leq 2 M S_R \leq M \left( 2 + \log_2 q \right)~.
\]
Thus, $\kappa \le (8+\log_2 q)M$.

Now, similarly to the above bound on $X_n^*$,
\[
\hX_n \geq - \kappa + \sum_{r=0}^R \sum_{s=0}^{S_r - S_{r-1}} \hX^{({r,s})}
\]
so that the regret $\hX_n - X_n^*$ is larger than
\[
\hX_n - X_n^* \geq -2\kappa + \sum_{r=0}^R \sum_{s=0}^{S_{r} - S_{r-1}}
\left( \hX^{({r,s})} - X_{k_{(r,s)}}^{(r,s)} \right)~.
\]
Now note that each time step $t$ (but the last one) of epoch $(r,s)$ satisfies
$M_t \le M^{(r)}$ and $\eta_{(r,s)} \le 1/2M^{(r)}$. Therefore, we can apply
Lemma~\ref{lm:refined} to $\hX^{({r,s})} - X_{k_{(r,s)}}^{(r,s)}$ for each epoch
$(r,s)$. This gives
\[
\hX_n - X_n^* \geq -2\kappa -
\sum_{r=0}^R \sum_{s=0}^{S_{r} - S_{r-1}}
\left( \frac{\ln N}{\eta_{(r,s)}} + \eta_{(r,s)} Q_{k_{(r,s)}}^{({r,s})} \right)~.
\]
By definition of the algorithm, for all epochs $(r,s)$,
\[
    Q_{k_{(r,s)}}^{(r,s)} \le Q_{k_{(r,s)},t_{(r,s)} - 1} = Q^*_{t_{(r,s)} - 1}
\le
    4^{S_{r-1}+s} \Bigl(M^{(r)}\Bigr)^2
\]
and
\[
    \eta_{(r,s)} \le \sqrt{\ln N}\Big/\Bigl(2^{S_{r-1}+s} M^{(r)}\Bigr)~.
\]
Therefore,
\begin{eqnarray}
\nonumber
\lefteqn{
    \sum_{r=0}^R \sum_{s=0}^{S_{r} - S_{r-1}} \eta_{(r,s)} Q_{k_{(r,s)}}^{({r,s})}
}
\\*
\nonumber
& \le &
    \sum_{r=0}^R \sum_{s=0}^{S_{r} - S_{r-1}} 2^{S_{r-1}+s} M^{(r)} \sqrt{\ln N}
\\
\nonumber
& \le &
    \sum_{r=0}^R \sum_{s=1}^{S_{r} - S_{r-1}} 2^{S_{r-1}+s} (2M) \sqrt{\ln N}
    + \sum_{r=0}^R 2^{S_{r-1}} M^{(r)} \sqrt{\ln N}
\\
\nonumber
& \le &
    (2M) \sum_{s=1}^{S_R} 2^{s} \sqrt{\ln N} + 2^{S_R} \sum_{r=0}^R M^{(r)} \sqrt{\ln N}
\\
\label{SR=0-1}
& \le &
    (2M) 2^{S_R+1} \sqrt{\ln N} + 2^{S_R} (4 M) \sqrt{\ln N}
\\ \nonumber
&& \quad
    \mbox{\rm (using~(\ref{sum2power}) and $M^{(R)} \le 2M$)}
\\*
\nonumber
& \le &
    (16M) \sqrt{q\ln N}
\end{eqnarray}
since $q \ge 4^{S_R-1}$ implies $2^{S_R} \le 2\sqrt{q}$.

We now turn our attention to the remaining sum
\[
    \sum_{r=0}^R \sum_{s=0}^{S_{r} - S_{r-1}} \frac{\ln N}{\eta_{(r,s)}}~.
\]
By definition of the algorithm,
\[
    \eta_{(r,s)} = \left\{ \begin{array}{cl}
    \displaystyle{1/(2M^{(r)})}
    & \mbox{\rm if $S_{r-1} + s \leq \lceil (\log_2 \ln N)/2 \rceil$}
\\
    \displaystyle{\sqrt{\ln N}/ \left( 2^{S_{r-1}+s} M^{(r)} \right)}
    & \mbox{\rm otherwise.}
    \end{array} \right.
\]
We denote by $(r^*,s^*)$ the last couple $(r,s)$ for which
$\eta_{r,s} = 1/(2M^{(r)})$. With obvious notation, a crude
overapproximation leads to
\begin{eqnarray}
\nonumber
\lefteqn{\sum_{r=0}^R \sum_{s=0}^{S_{r} - S_{r-1}} \frac{\ln N}{\eta_{(r,s)}}} \\
\nonumber
& & \leq \sum_{(r,s) \leq (r^*,s^*)} 2 M^{(r)} \ln N
+ \sum_{r=0}^R \sum_{s=0}^{S_{r} - S_{r-1}} 2^{S_{r-1}+s} M^{(r)} \sqrt{\ln N}~.
\end{eqnarray}
We already have the upper bound $(16M)\sqrt{q\ln N}$ for the second sum.
For the first one, we write
\begin{eqnarray*}
\lefteqn{
    \sum_{(r,s) \leq (r^*,s^*)} 2 M^{(r)} \ln N
}
\\ & = & \sum_{r=0}^{r^*} 2 M^{(r)} \ln N + \sum_{r=0}^{r^*-1} (S_r - S_{r-1}) \left( 2 M^{(r)} \right) \ln N \\
& & +               s^*  \left( 2 M^{(r^*)} \right) \ln N \\
& \leq & \sum_{r=0}^{R} 2 M^{(r)} \ln N + (S_{r^*-1} + s^*) (4M) \ln N \\
& \leq & 2M (\ln N ) \left( 3 + 2 \lceil (\log_2 \ln N)/2 \rceil \right)
\end{eqnarray*}
where we used~(\ref{sumR}).
The proof is concluded in the case $S_R \geq 1$ by putting things together
and performing some overapproximation.

% To complete the 'otherwise' part above.
% This also explains the max {1, ...} in the bound.
When $S_R = 0$, $q = 1$, $\kappa$ is simply less than $6M$,
(\ref{SR=0-1}) is less than $8M\sqrt{\ln N}$,
so that the bound holds as well in this case.
\qed
\end{pf}

\section*{Proof of Lemma~\ref{LmEWA_incr}}
We first note that Jensen's inequality implies that $\Phi$ is nonnegative.

The proof below is a simple modification of an argument first proposed in~\inlinecite{AuCeGe02}.
Note that we consider real-valued (non necessarily nonnegative) payoffs in what follows.
For $t = 1,\ldots,n$, we rewrite $p_{i,t} = w_{i,t}/W_t$, where
$w_{i,t} = e^{\eta_t X_{i,t-1}}$
and
$W_t = \sum_{j = 1}^{N} w_{j,t}$ (the payoffs $X_{i,0}$
are understood to equal 0, and thus, $\eta_1$ may be any positive number
satisfying $\eta_1 \geq \eta_2$).
Use
$w_{i,t}' = e^{\eta_{t-1} X_{i,t-1}}$ to denote the weight
$w_{i,t}$ where the parameter $\eta_{t}$ is replaced by $\eta_{t-1}$.
The associated normalization factor will be denoted by
$W_t' = \sum_{j=1}^N w_{j,t}'$.
Finally, we use $j_t^*$ to denote the expert with the largest cumulative
payoff after the first
$t$ rounds (ties are broken by choosing the expert with
smallest index). That is, $X_{j_t^*,t} = \max_{i \le N} X_{i,t}$.
We also make use of the following technical lemma.
\begin{lemma}[\opencite{AuCeGe02}]
\label{tech-varying}
For all $N \ge 2$, for all $\beta \ge \alpha \ge 0$, and for all
$d_1,\ldots,d_N \ge 0$ such that $\sum_{i=1}^N e^{-\alpha d_i} \ge 1$,
\[
    \ln\frac{\sum_{i=1}^N e^{-\alpha d_i}}{\sum_{j=1}^N e^{-\beta d_j}}
\le
    \frac{\beta-\alpha}{\alpha}\ln N~.
\]
\end{lemma}
\begin{pf*}{Proof (of Lemma \ref{tech-varying})}
We begin by writing
\begin{eqnarray*}
    \ln\frac{\sum_{i=1}^N e^{-\alpha d_i}}{\sum_{j=1}^N e^{-\beta d_j}}
&=&
    \ln\frac{\sum_{i=1}^N e^{-\alpha d_i}}{\sum_{j=1}^N e^{(\alpha-\beta)d_j} e^{-\alpha d_j}}
\\ &=&
    - \ln \esp{}{e^{(\alpha-\beta) D}}
\\ & \le &
    (\beta-\alpha)\esp{}{D}
\end{eqnarray*}
where we applied Jensen inequality to the random variable $D$ taking value $d_i$ with probability
$e^{-\alpha d_i}/\sum_{j=1}^N e^{-\alpha d_j}$ for each $j=1,\ldots,N$.
Since $D$ takes at most $N$ distinct values,
its entropy $H(D)$ is at most $\ln N$. Therefore
\begin{eqnarray*}
    \ln N \ge H(D)
&=&
    \frac{\sum_{i=1}^N e^{-\alpha d_i}}{\sum_{j=1}^N e^{-\beta d_j}}
    \left( \alpha d_i + \ln \sum_{j=1}^N e^{-\beta d_j} \right)
\\* &=&
    \alpha\esp{}{D} + \ln \sum_{j=1}^N e^{-\beta d_j}
\ge
    \alpha\esp{}{D}
\end{eqnarray*}
where the last inequality holds since $\sum_{i=1}^N e^{-\alpha d_i} \ge 1$.
Hence $\esp{}{D} \le (\ln N)/\alpha$. As $\beta > \alpha$ by hypothesis,
we can plug
the bound on $\esp{}{D}$ in the upper bound
above and conclude the proof.
\qed
\end{pf*}

\begin{pf*}{Proof of Lemma \ref{LmEWA_incr}}
As it is  usual in the analysis of the  exponentially
weighted average predictor, we study
the evolution of $\ln(W_{t+1}/W_t)$. However, here
we need to couple this term with $\ln(w_{j_{t-1}^*,t}/w_{j_t^*,{t+1}})$
including in both terms the time-varying parameters $\eta_t, \, \eta_{t+1}$.
Tracking the currently best expert $j_t^*$ is used to lower
bound the weight $\ln(w_{j_t^*,{t+1}}/W_{t+1})$. In fact,
the weight of the overall best expert (after $n$ rounds) could
get arbitrarily small during the prediction process.
We thus obtain the following
\begin{eqnarray*}
\lefteqn{
    \frac{1}{\eta_t}\ln\frac{w_{j_{t-1}^*,t}}{W_t}
    - \frac{1}{\eta_{t+1}}\ln\frac{w_{j_t^*,{t+1}}}{W_{t+1}}
}
\\* &=&
    \left(\frac{1}{\eta_{t+1}} - \frac{1}{\eta_t}\right)\ln\frac{W_{t+1}}{w_{j_t^*,{t+1}}}
+
    \frac{1}{\eta_t}\ln\frac{w_{j_t^*,{t+1}}'/W_{t+1}'}{w_{j_t^*,{t+1}}/W_{t+1}}
+
    \frac{1}{\eta_t}\ln\frac{w_{j_{t-1}^*,t}/W_t}{w_{j_t^*,{t+1}}'/W_{t+1}'}
\\* &=& (A) + (B) + (C)~.
\end{eqnarray*}
We now bound separately the three terms on the right-hand side.
The term~$(A)$ is easily bounded
by using
$\eta_{t+1} \leq \eta_t$ and using the fact that $j_t^*$ is the index
of the expert with largest payoff after the first $t$ rounds. Therefore,
$w_{j_t^*,{t+1}}/W_{t+1}$ must be at least $1/N$. Thus we have
\[
    (A) = \left(\frac{1}{\eta_{t+1}} - \frac{1}{\eta_t}\right)\ln\frac{W_{t+1}}{w_{j_t^*,{t+1}}}
\le
    \left(\frac{1}{\eta_{t+1}} - \frac{1}{\eta_t}\right)\ln N~.
\]
We proceed to bounding the term~$(B)$ as follows
\begin{eqnarray*}
(B) &=& \frac{1}{\eta_t}\ln\frac{w_{j_t^*,{t+1}}'/W_{t+1}'}{w_{j_t^*,{t+1}}/W_{t+1}}
=
\frac{1}{\eta_t}\ln\frac{\sum_{i=1}^N e^{- \eta_{t+1} (X_{j_{t}^*,t}-X_{i,t})}}{
    \sum_{j=1}^N e^{-\eta_t (X_{j_{t}^*,t}-X_{j,t})} }
\\ & \le &
    \frac{\eta_t-\eta_{t+1}}{\eta_t\eta_{t+1}} \ln N
=
    \left(\frac{1}{\eta_{t+1}}-\frac{1}{\eta_t}\right)\ln N
\end{eqnarray*}
where the inequality is proven by applying Lemma~\ref{tech-varying} with
$d_i = X_{j_{t}^*,t} - X_{i,t}$.
Note that $d_i \ge 0$ since $j_t^*$ is the index of
the expert with largest payoff after the first $t$ rounds and
$\sum_{i=1}^N e^{-\eta_{t+1} d_i} \ge 1$
as for $i=j_{t}^*$ we have $d_i=0$.  \\
The term~$(C)$
is first split as follows,
\[
    (C) = \frac{1}{\eta_t}\ln\frac{w_{j_{t-1}^*,t}/W_t}{w_{j_t^*,{t+1}}'/W_{t+1}'}
=
    \frac{1}{\eta_t} \ln\frac{w_{j_{t-1}^*,t}}{w_{j_t^*,{t+1}}'}
    + \frac{1}{\eta_t}\ln\frac{W_{t+1}'}{W_t}~.
\]
We bound separately each one of the two terms on the right-hand side.
For the first one, we have
\[
    \frac{1}{\eta_t} \ln\frac{w_{j_{t-1}^*,t}}{w_{j_t^*,{t+1}}'}
=
    \frac{1}{\eta_t}
    \ln\frac{e^{\eta_t X_{j_{t-1}^*,t-1}}}{e^{\eta_t X_{j_t^*,t}}}
=
    X_{j_{t-1}^*,t-1} - X_{j_t^*,t}~.
\]
The second term is handled by using the very definition of
$\Phi$,
\begin{eqnarray*}
    \frac{1}{\eta_t}\ln\frac{W_{t+1}'}{W_t}
=
    \frac{1}{\eta_t}
    \ln\frac{\sum_{i=1}^N w_{i,t} e^{\eta_t x_{i,t}}}{W_t}
& = & \frac{1}{\eta_t}
    \ln {\sum_{i=1}^N p_{i,t} e^{\eta_t x_{i,t}}}
\\ & = &
    \sum_{i=1}^N p_{i,t} x_{i,t} + \Phi(\bp_t, \, \eta_t, \, \bx_t)~.
\end{eqnarray*}
Finally, we plug back in the main equation the bounds on the first two
terms $(A)$ and $(B)$, and the bounds on the two parts of the term $(C)$.
After rearranging we obtain
\begin{eqnarray*}
    0
& \le &
    \left( X_{j_{t-1}^*,t-1} - X_{j_t^*,t} \right) + \sum_{i=1}^N p_{i,t} x_{i,t}
+ \Phi(\bp_t, \, \eta_t, \, \bx_t)
\\ && \quad - \;
    \frac{1}{\eta_{t+1}}\ln\frac{w_{j_t^*,{t+1}}}{W_{t+1}}
    + \frac{1}{\eta_t}\ln\frac{w_{j_{t-1}^*,t}}{W_t}
\\ && \quad + \;
    2\left(\frac{1}{\eta_{t+1}}-\frac{1}{\eta_t}\right)\ln N~.
\end{eqnarray*}
We apply the above inequalities to each $t=1,\ldots,n$ and sum up using
\begin{eqnarray*}
 && \sum_{t=1}^n X_{j_{t-1}^*,t-1} - X_{j_t^*,t}
=
- \max_{j = 1,\ldots,N} X_{j,n} \\
  \mbox{and}  && \sum_{t=1}^n \left( - \frac{1}{\eta_{t+1}}\ln\frac{w_{j_t^*,{t+1}}}{W_{t+1}}
    + \frac{1}{\eta_t}\ln\frac{w_{j_{t-1}^*,t}}{W_t} \right)
\le
    - \frac{1}{\eta_1}\ln\frac{w_{j_0^*,1}}{W_1} = \frac{\ln N}{\eta_1}
\end{eqnarray*}
to conclude the proof.
\qed
\end{pf*}

\end{article}

\begin{thebibliography}{}

% Erased references correspond to the bandit problems

%\bibitem[\protect\citeauthoryear{Allenberg-Neeman and Auer}{2004}]{AA04}
%C.~Allenberg-Neeman and P.~Auer.
%\newblock Personal communication.

\bibitem[\protect\citeauthoryear{Allenberg-Neeman and Neeman}{2004}]{AlNe04}
C.~Allenberg-Neeman and B.~Neeman.
\newblock Full information game with gains and losses.
\newblock Algorithmic Learning Theory, 15th International Conference, ALT 2004, Padova,
Italy, October 2004, Proceedings, volume 3244 of Lecture Notes in Artificial Intelligence, pages 264-278. Springer, 2004.

\bibitem[\protect\citeauthoryear{Auer et al.}{2002}]{AuCeFrSc02}
P.~Auer, N.~Cesa-Bianchi, Y.~Freund, and R.E.~Schapire.
\newblock The nonstochastic multiarmed bandit problem.
\newblock \emph{SIAM Journal on Computing}, 32:48--77, 2002.

%GS Should be Auer et al. according to Machine Learning standards--but we already have an Auer et al. (2002)
\bibitem[\protect\citeauthoryear{Auer, Cesa-Bianchi, and Gentile}{2002}]{AuCeGe02}
P.~Auer, N.~Cesa-Bianchi, and C.~Gentile.
\newblock Adaptive and self-confident on-line learning algorithms.
\newblock \emph{Journal of Computer and System Sciences}, 64:48--75, 2002.

\bibitem[\protect\citeauthoryear{Cesa-Bianchi et al.}{1997}]{CFHHSW97}
N.~Cesa-Bianchi, Y.~Freund, D.P.~Helmbold, D.~Haussler, R.~Schapire, and M.K.~Warmuth.
\newblock How to use expert advice.
\newblock \emph{Journal of the ACM}, 3:427--485, 1997.

\bibitem[\protect\citeauthoryear{Cesa-Bianchi and Lugosi}{2003}]{CeLu03}
N.~Cesa-Bianchi and G.~Lugosi.
\newblock Potential-based algorithms in on-line prediction and game theory.
\newblock \emph{Machine Learning}, 51:239--261, 2003.

\bibitem[\protect\citeauthoryear{Cesa-Bianchi and Lugosi}{2006}]{CBL05}
N.~Cesa-Bianchi and G.~Lugosi.
\newblock {\em Prediction, Learning, and Games}.
\newblock Cambridge University Press, 2006.

\bibitem[\protect\citeauthoryear{Cesa-Bianchi, Lugosi, and Stoltz}{2005}]{CeLuSt04}
N.~Cesa-Bianchi, G.~Lugosi, and G.~Stoltz.
\newblock Minimizing regret with label efficient prediction.
\newblock {\em IEEE Transactions on Information Theory}, 51:2152--2162, 2005.

\bibitem[\protect\citeauthoryear{Cesa-Bianchi, Lugosi, and Stoltz}{2004}]{CLS04}
N.~Cesa-Bianchi, G.~Lugosi, and G.~Stoltz.
\newblock Regret minimization under partial monitoring.
\newblock Submitted for journal publication, 2004.

%GS Not needed anymore, I think
%\bibitem[\protect\citeauthoryear{Cesa-Bianchi, Mansour, and Stoltz}{2005}]{CeMaSt05}
%N.~Cesa-Bianchi, Y.~Mansour, and G.~Stoltz.
%\newblock Improved second-order bounds for prediction with expert advice.
%\newblock In {\em Proceedings of the 18th Annual Conference on Computational
%Learning Theory}, pages 217--232, 2005.

\bibitem[\protect\citeauthoryear{Freedman}{1975}]{Fre75}
D.~A.~Freedman.
\newblock On tail probabilities for martingales.
\newblock {\em The Annals of Probability}, 3:100--118, 1975.

\bibitem[\protect\citeauthoryear{Freund and Schapire}{1997}]{FS97}
Y.~Freund and R.E.~Schapire.
\newblock  A decision-theoretic generalization of on-line learning and an application to boosting.
\newblock {\em Journal of Computer and System Sciences}, 55(1):119--139, 1997.

\bibitem[\protect\citeauthoryear{Helmbold et al.}{1998}]{HeScSiWa98}
D.~P. Helmbold, R.~E. Schapire, Y.~Singer, and M.~K. Warmuth.
\newblock On-line portfolio selection using multiplicative updates.
\newblock {\em Mathematical Finance}, 8:325--344, 1998.

%\bibitem[\protect\citeauthoryear{Hart and Mas-Colell}{2001}]{HM}
%S. Hart and A. Mas-Colell.
%\newblock A Reinforcement Procedure Leading to Correlated Equilibrium.
%\newblock Economic Essays, Gerard Debreu, Wilhelm Neuefeind, and Walter Trockel (editors), 181-200. Springer, 2001.

\bibitem[\protect\citeauthoryear{Littlestone and Warmuth}{1994}]{LW94}
N.~Littlestone and M.K.~Warmuth.
\newblock The weighted majority algorithm.
\newblock {\em Information and Computation}, 108:212--261, 1994.

\bibitem[\protect\citeauthoryear{Piccolboni and Schindelhauer}{2001}]{PiSc01}
A.~Piccolboni and C.~Schindelhauer.
\newblock Discrete prediction games with arbitrary feedback and loss.
\newblock In {\em Proceedings of the 14th Annual Conference on Computational
Learning Theory}, pages 208--223, 2001.

\bibitem[\protect\citeauthoryear{Vovk}{1998}]{Vov98}
V.G.~Vovk.
\newblock A game of prediction with expert advice.
\newblock {\em Journal of Computer and System Sciences}, 56(2):153--73, 1998.

\end{thebibliography}
\end{document}